\input amstex
\documentstyle{amsppt}
\pageno=1
\magnification1200
\catcode`\@=11
\def\logo@{}
\catcode`\@=\active
\NoBlackBoxes
\vsize=23.5truecm
\hsize=16.5truecm

\def\d{d\!@!@!@!@!@!{}^{@!@!\text{\rm--}}\!}

\def\crp{\overline{\Bbb R}_+}

\def\rnp{{\Bbb R}^n_+}
\def\rnm{\Bbb R^n_-}

\def\crnp{\overline{\Bbb R}^n_+}

\def\ang#1{\langle {#1} \rangle}

\def\ttilde{\overset{\,\approx}\to}

\def\Zfrac{\tsize\frac1{\raise 1pt\hbox{$\scriptstyle z$}}}

\def\rp{ \Bbb R_+}

\define\tr{\operatorname{tr}}
\define\op{\operatorname{OP}}
\define\opg{\operatorname{OPG}}
\define\opt{\operatorname{OPT}}
\define\opk{\operatorname{OPK}}
\define\Tr{\operatorname{Tr}}
\define\TR{\operatorname{TR}}

\define\stimes{\!\times\!}
\def\res{\operatorname {res}}

\def\slint{{\int\!\!\!\!\!\!@!@!@!@!{-}}}
\def\tslint{{\int\!\!\!\!\!@!@!@!{-}}}

\document

\topmatter
\title
{Traces and quasi-traces on the Boutet de Monvel algebra}
\endtitle
\author{Gerd Grubb and Elmar Schrohe}
\endauthor
\address
{Copenhagen University, Mathematics Department,
Universitetsparken 5, 2100 Copenhagen, Denmark.
E-mail {\tt grubb\@math.ku.dk}}
\endaddress
\address
{Institut f\"ur Mathematik, Universit\"at Hannover,  Welfengarten 1,
30167 Hannover, Germany.}
E-mail {\tt schrohe\@math.uni-hannover.de}
\endaddress

\rightheadtext{Quasi-traces}

\abstract 
We construct an analogue of Kontsevich and Vishik's canonical trace
for pseudodifferential boundary value problems in the Boutet
de Monvel calculus on compact manifolds with boundary.
For an operator $A$ in the calculus (of class zero), and an auxiliary
operator $B$, formed of the Dirichlet realization of a strongly
elliptic second-order differential operator and an elliptic operator
on the boundary, we consider the coefficient $C_0(A,B)$ of
$(-\lambda)^{-N}$ in the asymptotic expansion of the resolvent trace
$\operatorname{Tr}(A(B-\lambda)^{-N})$ (with $N$ large) in powers and
log-powers of $\lambda$. This coefficient identifies with the zero-power
coefficient in the Laurent series for the zeta function 
$\operatorname{Tr}(AB^{-s})$ at $s=0$, when $B$ is invertible.
We show that $C_0(A,B)$ is in general a quasi-trace, in the sense that
it vanishes on commutators $[A,A']$ modulo local terms, and has a
specific value independent of the auxiliary operator, modulo local terms.
The local ``errors'' vanish when $A$ is a singular Green operator of 
noninteger order, or of integer order with a certain parity; then
$C_0(A,B)$ is a trace of $A$. They do not in general vanish when the 
interior ps.d.o. part of $A$ is nontrivial.


\noindent{\it Math.~classification:} 58J42 -- 35S15
\endabstract
\endtopmatter

\head  1. Introduction \endhead

In their work on the geometry of determinants of elliptic operators, 
Kontsevich and Vishik 
\cite{KV} introduced the
{\it canonical trace},
a novel functional $\operatorname{TR} A$ defined for a class of
pseudodifferential  operators ($\psi$do's) $A$  
on a closed $n$-dimensional manifold $X$. 
What makes it remarkable is that  
it extends the standard operator trace, and hereby complements 
the noncommutative residue. 

The  noncommutative residue $\operatorname{res}A$ was discovered 
by Wodzicki \cite{W} and, independently, by Guillemin \cite{Gu};
it is a trace (i.e., a nontrivial linear functional which vanishes 
on commutators) on the full algebra ${\Cal A}$
of all classical pseudodifferential operators on $X$.
The noncommutative residue turns out to be the only 
trace on $\Cal A$ with this property, up to multiples.
Moreover, the value of $\operatorname{res}A$ is determined from finitely many 
terms in the asymptotic expansion of  the symbol of $A$ 
(in fact, from the component of homogeneity  $-n$ in this expansion); 
we call such functionals ``locally determined''.

It vanishes on operators of noninteger order 
or of order $<-n$, and it is here that Kontsevich and Vishik's
functional takes over: $\TR A$ is well-defined and nontrivial when $A$ is of
noninteger order or of order $<-n$, 
and it vanishes on commutators $[A,A']$
of these types. It is ``global''  in that the value of $\TR A$ 
depends on the full operator and cannot be determined from the terms
in the asymptotic expansion of the symbol. The
canonical trace equals the standard trace when the order is $<-n$; 
it is moreover defined on integer-order operators with certain parity
properties, see below.  

In this article we shall generalize the canonical
trace to pseudodifferential boundary value problems 
in the calculus of Boutet de Monvel \cite{B}.

The original approach of \cite{KV} is based on studies of
generalized zeta func\-tions \linebreak$\zeta (A,P,s)=\Tr(AP^{-s})$ and certain
regularizations of them, where $P$ is an auxiliary elliptic operator.
This does not extend easily to
the situation of manifolds with boundary;
in the case without boundary, $P^{-s}$ is a classical $\psi $do again
(by Seeley \cite{S}),
whereas complex powers of elliptic
boundary problems fall outside the calculus of Boutet de Monvel.
However, the introduction of $\TR A$ can instead be based on
trace expansions of heat operators $ Ae^{-tP}$ (Lesch \cite{L}) or resolvents
$A(P-\lambda )^{-N}$ (Grubb \cite{G4}), where the latter
admit a direct generalization to manifolds with boundary (Grubb and
Schrohe \cite{GSc}).

Therefore, let us explain the functional of \cite{KV} from the
resolvent point of view:
Let $A$ be a classical $\psi $do of order
$\nu$; its symbol has an expansion in local coordinates $a\sim
\sum_{j\ge 0}a_{\nu -j}$
with $a_{\nu -j}(x,\xi )$ smooth in $(x,\xi )$ and homogeneous of
degree
$\nu -j$ in $\xi $ for $|\xi |\ge 1$. 
When $P$ is an auxiliary elliptic
operator of integer order
$m>0$ with no
principal symbol eigenvalues on $\Bbb R_-$, 
a calculation in local coordinates shows that
the operator family
$A(P-\lambda )^{-N}$ has for $N>(n+\nu )/m$ a trace expansion $$
\Tr (A(P-\lambda )^{-N})
\sim
\sum_{j\ge 0}\tilde c_{j}(-\lambda )^{\frac {n+\nu -j}m-N}
+\sum_{k\ge 0}
\bigl(\tilde c'_{k}\log (-\lambda )+\tilde
c''_{k}\bigr)(-\lambda )^{-k-N},\tag 1.1$$
for
$\lambda \to \infty$ in a small sector around $\Bbb R_-$,
cf\. Grubb and Seeley \cite{GS1}. 
Here each $\tilde c_j$ (and each $\tilde c_k'$) comes from a specific
homogeneous term in the symbol of $A(P-\lambda)^{-N}$, whereas the 
$\tilde c''_k$ depend on the full symbol. So the
coefficients $\tilde c_j$ and $\tilde c_k'$ depend each on a finite
set of homogeneous terms in the symbols of $A$ and $P$; we call such
coefficients `locally determined' (or `local'), while the $\tilde
c''_k$ are called `global'. When $\nu \notin
\Bbb Z$, the $\tilde c'_k$ vanish. When $\nu \in\Bbb Z$ and $(j-n-\nu
)/m$ is an integer $k\ge 0$, both $\tilde c_j$ and $\tilde c''_k$
contribute to the power $(-\lambda )^{-k-N}$; their sum is
independent of the choice of local coordinates, whereas the
splitting in $\tilde c_j$ and $\tilde c''_k$ depends in a well-defined way
on the symbol structure in local coordinates (see  \cite{GS1,
Th.\ 2.1} or \cite{G4, Th. 1.3}).

In (1.1), 
$m\cdot \tilde c_0'= \operatorname{res}A$, depending solely on $A$,
cf.\ \cite{W}, \cite{Gu}. Moreover, cf.\ \cite{KV}, \cite{G4},
$$
\operatorname{TR} A= \tilde c''_0\tag1.2
$$
in the following four cases:
\roster
\item "(1)" $\nu <-n$,
\item "(2)" $\nu \notin\Bbb Z$,
\item "(3)" $\nu \in\Bbb Z$, $A$ is even-even and $n$ is odd,
\item "(4)" $\nu \in\Bbb Z$, $A$ is even-odd and $n$ is even;
\endroster
in the cases (3) and (4), $P$ is taken to be even-even with $m$ even.
In all the cases
(1)--(4), $\operatorname{res}A=0$, there is no contribution to 
$(-\lambda)^{-N}$ from the first sum in (1.1),  and
$$
\tilde c''_0=\int\slint \tr a(x,\xi )\,\d\xi dx,\tag1.3
$$
where $\tslint \tr a(x,\xi )\,\d\xi $ is a finite part integral (explained
in detail in Section 3 below).

When $\nu \in\Bbb Z $, we say that
$A$ or $a$ has {\it even-even}
alternating
parity (in short: is even-even), when the
symbols with even (resp\. odd) degree $\nu -j$ are even (resp\. odd)
in $\xi $:
$$
a_{\nu -j}(x,-\xi )=(-1)^{\nu -j}a_{\nu -j}(x,\xi ) \text{ for }|\xi
|\ge 1,\tag1.4
$$
and the derivatives in $x$ and $\xi $ likewise have this property.
$A$ and $a$ are said to have {\it even-odd} alternating
parity in the reversed
situation, where the
symbols with even (resp\. odd) degree $\nu -j$ are odd (resp\. even)
in $\xi $:
$$
a_{\nu -j}(x,-\xi )=(-1)^{\nu -j-1}a_{\nu -j}(x,\xi ) \text{ for
}|\xi |\ge 1,\tag1.5
$$
etc.
These properties are preserved under coordinate changes.
For brevity, we shall say that $A$ (or $a$) {\it has a parity that
fits with the
dimension $n$}, when (3) or (4) holds.

\cite{KV} only considered the cases (1)--(3); the operators
satisfying (3) were called odd-class operators. (4) was
included in \cite{G4}.

The relations can also be formulated in terms of the generalized
zeta function
\linebreak$\zeta (A,P,s)=\Tr(AP^{-s})$; here we assume $P$ invertible for
simplicity. It is known from \cite{W}, \cite{Gu} that $\zeta
(A,P,s)$ extends meromorphically from large $\operatorname{Re}s$
across $\operatorname{Re}s=0$ with a simple pole at $0$,
$$
\zeta (A,P,s)=C_{-1}(A,P)s^{-1}+C_0(A,P)+O(s) \text{ for }s\to 0;\tag1.6
$$
this follows also by translating (1.1) to a statement on the pole
structure of $\zeta (A,P,s)$ (as e.g.\ in \cite{GS2, Cor\. 2.10}).
Here$$\aligned
C_{-1}(A,P)&=\tilde c'_0=\tfrac 1m \operatorname{res}A,\\
C_0(A,P)&=\tilde c_{\nu +n}+\tilde c''_0,
\endaligned\tag1.7$$
with $\tilde c_{\nu +n}$ defined as 0 if $\nu <-n$ or $\nu
\notin\Bbb Z$. In the cases (1)--(4), $\tilde c_{\nu +n}$ and
$\operatorname{res}A$ vanish 
(for any choice of local coordinates),
so that$$
\operatorname{TR}A=
\tilde c''_0 =C_0(A,P)
=\zeta (A,P,0).\tag1.8
$$

The coefficient $C_0(A,P)$ has an interest also when none of the
conditions (1)--(4) is satisfied. Then it has the properties:$$
C_0(A,P)-C_0(A,P') \text{ and }C_0([A,A'],P)\text{ are locally
determined}.\tag1.9
$$
Moreover, $C_0(A,P)$ can be written as (1.3) plus
local terms.
Functionals 
with properties as in (1.9) for a system of auxiliary elliptic operators
$P$ will be called quasi-traces.

By use of more functional
calculus, defining $\log P$, one can moreover show that the two
expressions in (1.9) are noncommutative residues of suitable
combinations of the given operators and $\log P$; cf\. \cite{KV}
and Okikiolu \cite{O} for $C_0(A,P)-C_0(A,P')$, cf.\ Melrose and
Nistor \cite{MN}
for $C_0([A,A'],P)$. (In works of
Melrose et al., notation such as $\widehat{\operatorname{Tr}}(A)$ or
$\Tr_P(A)$ is used for $C_0(A,P)$; it is called a regularized
trace there.)

\medskip
In order to extend the definition of
the canonical trace to operators in the calculus of
Boutet de Monvel \cite{B} on an $n$-dimensional compact
$C^\infty$-manifold $X$ with boundary $\partial X=X'$,
we shall rely on resolvent expansions analogous to
(1.1). To this end we choose an auxiliary operator
$P_{1,\operatorname{D}}$,
which is the (invertible) Dirichlet realization of
a strongly elliptic principally scalar second-order differential operator.

It was shown in \cite{GSc} that when $A=P_++G$ is a pseudodifferential
boundary operator ($\psi$dbo) in this calculus of order
$\nu \in\Bbb Z$ with $G$ of class 0, then
there is a trace expansion for $N>(n+\nu )/2$:
$$
\Tr(A(P_{1,\operatorname{D}}-\lambda )^{-N})\sim \sum_{j\ge
0}\!\tilde c_{j}(-\lambda)
^{\frac{n+\nu - j}{2}-N}
+\sum_{k\ge
0}(\tilde c'_k\log (-\lambda )+\tilde c''_k)(-\lambda )
^{-\frac k2-N}.
\tag 1.10$$
The proof depends on a reduction to parameter-dependent $\psi $do's
on $X'$, where \cite{GS1} can be used; it is easily generalized to
cases where $A=G$ of noninteger order $\nu $ and class 0. When $\nu
\notin\Bbb Z$, the $\tilde c'_k$ vanish.

As explained e.g.\ in \cite{GS2} this yields the expansions:
$$\align
\Tr(Ae^{-tP_{1,\operatorname{D}}})\sim
\sum_{j\ge 0}
c_{j} t
^{\frac{j-n-\nu  }{2}}
+\sum_{k\ge 0}(-  c'_k\log t +  c''_k)t
^{\frac k2},\tag1.11\\
 \Gamma (s)\Tr
(AP_{1,\operatorname{D}}^{-s})\sim \sum_{j\ge 0}\frac{
c_j }
{s+\frac{j-n-\nu  }{2}}
+\sum_{k\ge 0}\Bigl(\frac{  c'_k}{(s+\frac k2)^2} +\frac{
c''_k}
{s+\frac k2}\Bigr) .\tag1.12
\endalign
$$
In (1.10), $\lambda \to\infty $ in a sector of $\Bbb C$;
in (1.11), $t\to 0+$; (1.12) describes the pole structure of the
meromorphic extension of $\Gamma (s)\Tr
(AP_{1,\operatorname{D}}^{-s})$ to $\Bbb C$.
The coefficients $c_j,  c_k', c''_k$ are
proportional to the coefficients $\tilde c_j, \tilde c_k', \tilde
c''_k$ in (1.10) by
universal constants; in particular,$$
c'_0=\tilde c'_0,\quad c''_0=\tilde c''_0, \quad c_{\nu +n}=\tilde
c_{\nu +n},\tag1.13
$$
where $c_{\nu +n}$ and $\tilde c_{\nu +n}$ are taken as 0 when $\nu
<-n$ or $\nu \notin\Bbb Z$.
The coefficients $c_j$ and $c_k'$ are locally determined,
whereas the $c''_k$ are global.

Relation (1.12) implies that the
generalized zeta function $\zeta
(A,P_{1,\operatorname{D}},s)=\Tr (AP^{-s}_{1,\operatorname{D}})$
extends meromorphically with
$$\gathered
\zeta
(A,P_{1,\operatorname{D}},s)=C_{-1}(A,P_{1,\operatorname{D}})s^{-1}
+C_0(A,P_{1,\operatorname{D}})+O(s) \text{ for }s\to 0,\\
C_{-1}(A,P_{1,\operatorname{D}})=\tilde c'_0,\quad
C_0(A,P_{1,\operatorname{D}})=\tilde c_{\nu +n}+\tilde c''_0.
\endgathered\tag1.14
$$
In  \cite{FGLS},
an analogue of Wodzicki's noncommutative residue
was introduced by Fedosov, Golse, Leichtnam, and Schrohe
for operators in Boutet de Monvel's calculus.
In \cite{GSc} we were able to show that the coefficient $\tilde
c_0'=C_{-1}(A,P_{1,\operatorname D})$ 
in the above expansions satisfies the relation
$$\tilde c_0'=\tfrac12 \res A.\tag1.15$$

We shall presently investigate $C_0(A,P_{1,\operatorname{D}})$ as a
candidate for a canonical trace.
We establish quasi-trace properties as in (1.9), with
formulas for the
``value modulo local terms''.
Moreover, we extract some cases where
$C_0(A,P_{1,\operatorname{D}})$ is independent of
$P_{1,\operatorname{D}}$ and vanishes
on commutators $[A,A']$, so that it is a trace.
Here we show in particular that
$C_0([A,A'],P_{1,\operatorname{D}})=0$
whenever $A$ and $A'$ are singular Green operators
of class zero and orders $\nu$ and $\nu'$
with $\nu+\nu'<1-n$ or $\nu+\nu'\notin \Bbb Z$;
the same is true
when $\nu+\nu'$ is an integer
and a certain parity holds, provided we
narrow down slightly the class of admissible auxiliary operators
$P_1$.
When $A=P_++G$ and $A'=P'_++G'$ both have nontrivial pseudodifferential
part (and hence the orders are integer), however, one cannot
hope for more than the quasi-trace property, see Remark 4.2 for a discussion
of this point.

In contrast with the boundaryless case, the powers
$P_{1,\operatorname{D}}^s$
are far from
belonging to the $\psi $dbo calculus when $s\notin \Bbb Z$, so we do not
expect to find residue formulas as in \cite{KV}, \cite{O}, \cite{MN}.

The quasi-trace property has an interest in
itself, for situations where one has some control over the local
terms. Our formulas for $C_0(P_+,P_{1,\operatorname{D}})$ in Theorem 4.1
and for  $C_0(G,P_{1,\operatorname{D}})$ in Theorem 3.6
moreover show how the values are related to formulas
for the boundaryless manifolds $\widetilde X$ resp\. $X'$.

\head 2. Preliminaries \endhead

Let $X$ be an $n$-dimensional compact
$C^\infty $ manifold with boundary $\partial X=X'$, $X$ and $X'$
provided with $C^\infty $ vector bundles
$E$ and $E'$. We can assume that $X$ is smoothly imbedded in an
$n$-dimensional manifold $\widetilde X$ without boundary, provided
with a vector bundle $\widetilde E$ such that $E=\widetilde E|_X$.
We consider the algebra of (one-step) polyhomogeneous
operators in the calculus of Boutet de Monvel 
$$
\pmatrix P_++G & K\\ \quad&\quad \\T&S
\endpmatrix
\:
\matrix
C^\infty (E)&&C^\infty(E)\\
\oplus&\longrightarrow&\oplus\\
 C^\infty (E')&& C^\infty (E')
\endaligned.
 \tag2.1$$
Here  $P$ is a pseudodifferential operator satisfying the
transmission condition.
The subscript `$+$' indicates that we are taking the truncation of
$P$ to $X$,
i.e., the operator given by extending a function in $C^\infty (E)$ by
zero to a
function in $L^2(\widetilde X,\widetilde E)$, applying $P$, and
restricting the result
to $X^\circ$; the transmission property assures that this gives an
element of $C^\infty (E)$.
Moreover, $G$  is a singular Green operator (s.g.o.), $K$ a Poisson
operator,
$T$ a trace
operator, and $S$ a pseudodifferential operator on $X'$.
We assume that all are of order $\nu$. Details on the calculus
can e.g\. be found in \cite{G2}.

As a first observation we note that
for classical (one-step polyhomogeneous) $\psi $do's $P$ of order
$\nu $ on
$\widetilde X$, the calculus of Boutet de Monvel requires $\nu \in\Bbb Z$
in
order for the transmission condition to be satisfied at $\partial X$.
The definition of the operators $G$, $K$
and $T$ (and, of course, $S$), however, extends readily to noninteger $\nu
$, so
one can also ask for an extension of the canonical trace
to  operators  of the form
$$
\pmatrix G & K\\ T&S
\endpmatrix ,\tag2.2$$
with pseudodifferential part equal to zero and all elements of order
$\nu \in \Bbb R$.

We next fix  an auxiliary operator:
We let $P_1$ be a second-order strongly elliptic differential
operator in $\widetilde E$
with scalar
principal symbol, and denote by $P_{1,D}$ its
Dirichlet realization on $X$. By possibly shifting the operator by a
constant,
we can assume that $P_1$ and $P_{1,\operatorname{D}}$ have positive
lower bound, so that the resolvents $Q_\lambda =(P_1-\lambda )^{-1}$
and $R_\lambda =(P_{1,\operatorname{D}}-\lambda )^{-1}$ are defined
for $\lambda $ in a region$$
\Lambda =\{ \,\lambda \in\Bbb C\mid \arg\lambda \in [\tfrac \pi
2-\varepsilon ,\tfrac{3\pi }2+\varepsilon] \text{ or }|\lambda |\le
r(\varepsilon )\,\}\tag2.3
$$
with suitable $\varepsilon>0$ and $r(\varepsilon)>0$.
Moreover, we take a strongly elliptic second-order
pseudodifferential operator $S_1$ on $C^\infty(X',E')$
and set
$
B=\pmatrix P_{1,D}& 0\\ 0&S_1\endpmatrix.
$
Then
$$\aligned
&\Tr \Bigg(\pmatrix P_++G & K\\ T&S
\endpmatrix(B-\lambda)^{-N}\Bigg)\\
&\qquad\qquad= \Tr_X \left((P_++G)(P_{1,D}-\lambda)^{-N}\right)
+\Tr_{X'}\left(S(S_1-\lambda)^{-N}\right)\endaligned\tag2.4
$$
is well-defined for
$N>(n+\nu)/2$. The behavior of
$\Tr_{X'}\left(S(S_1-\lambda)^{-N}\right)$ is well-known from the
theory of $\psi $do's on  manifolds without boundary,
so we shall restrict the attention to
$$\Tr_X \left((P_++G)(P_{1,D}-\lambda)^{-N}\right)\text{ if $\nu
\in\Bbb Z$},\text{ and }~~ \Tr_X
\left(G(P_{1,D}-\lambda)^{-N}\right)\text{ if $\nu
\in\Bbb R\setminus \Bbb Z$}.\tag2.5
$$
We henceforth denote $P_++G=A$.

A natural candidate for an extension of the canonical trace
is the functional $A\mapsto C_0(A,P_{1,\operatorname D})$ which associates
with the $\psi$dbo $A=P_++G$ the coefficient of $(-\lambda )^{-N}$ in
(1.10), equal to the coefficient of $s^0$ in the Laurent expansion
of $\zeta(A,P_{1,\operatorname D},s)$ at $s=0$ (1.14).
We shall show that this is indeed a good choice at least in two cases:

(i) $A=G$ is a singular Green operator of noninteger order $\nu$ and
class zero
with vanishing pseudodifferential part;

(ii)  $A=P_++G$ is of integer order $\nu$ with a pseudodifferential
operator
$P$ of normal order $0$ and $G$ of class zero (we say that a $\psi
$do symbol has normal order $d$ when it is $O(\xi _n^{d})$ at
the boundary).

There are two conditions involved, namely on the class of the singular
Green part and
on the normal order of the pseudodifferential part in (ii). They are
both natural, and related:

If $G$ is effectively of class $r>0$, i.e., if it can be written in the
form
$G=G_0+\sum_{j=0}^{r-1} K_j\gamma_j$, where
$G_0$ is a singular Green operator of class zero,
$\gamma_ju=(\partial _n^ju)|_{X'}$, and the $K_j$ are Poisson
operators with $K_{r-1}\ne 0$,
then $G$ will not be bounded on $L^2(E)$,
much less of trace class, even if its order
is arbitrarily low, cf.\ e.g.\ \cite{G2, Sect.\ 2.8}.
As we are looking for a functional coinciding with the standard
trace for operators of sufficiently low order,
we shall exclude operators of positive class.
An expansion analogous to (1.10) still
exists, with the sum over $k\ge 0$ replaced by a sum
over $k\ge k_0$ with
a possibly negative starting index $k_0$.
In this case, however, the coefficient $c'_0$
will not in general coincide with the noncommutative residue,
cf\.  \cite{GSc, Remark 1.2}.

For operators $P$ having the transmission property, we find that
$C_0(P_+, P_{1,\operatorname{D}})$ does indeed have a value modulo
local terms that is independent of $P_1$ modulo local terms.
The commutator $[P_+,P_+']$ of two pseudodifferential operators $P$
and $P'$ truncated to $X$, however,
will be of the form $P_+''+G''$, where, in general,
the singular Green operator $G''$ will be of class $r>0$
unless both
$P$ and $P'$ are of normal order $\le 0$, cf\. \cite{G2, Section 2.6}.
Hence we can only hope to show the commutator property for
this class.

Recall (cf\. e.g\.
\cite{G2, Lemma 1.3.1}) that any $\psi $do $P$ having the transmission
property at $X'$ can be written as a sum $$
P=P^{(1)}+P^{(2)}+P^{(3)},\tag 2.6
$$ where
$P^{(1)}$ is a differential operator,
$P^{(2)}$ is a $\psi $do whose symbol has  normal order $-1$, and the
symbol of $P^{(3)}$ has, near $X'$, a large power of $x_n$ as a
factor to the left or right.

For a {\it differential} operator $P^{(1)}$, the expansion
in (1.10) (with $A=P^{(1)}$) is valid without the second series over
$k$; all terms are locally determined.
Therefore $C_{-1}(P^{(1)},P_{1,\operatorname{D}})$
vanishes, and $C_0(P^{(1)},P_{1,\operatorname{D}})$ is
local in this case, and so is
$C_0(P^{(1)},P_{1,\operatorname{D}})-C_0(P^{(1)},P_{2,\operatorname{D}})$
for
another choice of auxiliary operator $P_2$. Thus
$C_0(P^{(1)},P_{1,\operatorname{D}})$ satisfies the first part of the
requirement for being a quasi-trace (with value zero modulo local terms).
As mentioned above, differential operators of positive normal order
will be left out from the commutator
considerations. So will
operators of type $P^{(3)}$, since the structure is not preserved under
differentiation of symbols.

At this point we can state a part of the results we will show in this
paper:

\proclaim{Theorem 2.1}
Let $A=P_++G$ and $A'=P'_++G'$ be of orders $\nu $ resp\. $\nu
'$, with $P$ resp\. $P'$ vanishing if $\nu $ resp\. $\nu ' \notin \Bbb
Z$, and with $G$ and $G'$ of class $0$. Let $P_1$ and $P_2$ be
two auxiliary operators as explained above.
Then

{\rm (i)}
$$
C_{0}(A,P_{1,\operatorname{D}})-C_{0}(A,P_{2,\operatorname{D}})\text{
is locally determined};\tag2.7
$$
it vanishes if $\nu <-n$ or $\nu \notin \Bbb Z$, and otherwise depends
solely on the terms of the first
$\nu +n+1$
homogeneity degrees in the symbols of $P_++G $, $P_1$ and $P_2$.

{\rm (ii)} If $\nu $ and $\nu '\in\Bbb Z$, assume that $P$ and $P'$
have normal order $0$; if $\nu $ or $\nu '\in \Bbb R\setminus \Bbb Z
$, assume that $P$ and $P'$ are zero. Then
$$
C_{0}([A,A'],P_{1,\operatorname{D}})\text{
is locally determined};\tag2.8
$$
it vanishes if $\nu +\nu '<-n$ or $\nu +\nu '\notin \Bbb Z$, and
otherwise depends
solely on the terms of the first
$\nu +\nu '+n+1$
homogeneity degrees in the symbols of $P_++G $, $P'_++G'$ and $P_1$.
\endproclaim

(Actually, the number of homogeneous terms entering from $G$ is one
step lower; see the statements in Section 3.)

The operators of the form $A=P_++G$ of integer order with $G$ of class
zero and $P$ of normal order $\le 0$ form an algebra;
Theorem 2.1 shows that $C_0(A,P_{1,\operatorname{D}})$ is a quasi-trace
on this algebra as well as on the singular Green operators of
class zero.

Similarly as in \cite{GSc}, our analysis is based on precise
information about the resolvent.
We therefore recall the structure of
$R_\lambda  =(P_{1,\operatorname{D}}-\lambda)^{-1}$:
The Dirichlet problem $$
(P_1-\lambda )u=f\text{ on }X,\quad \gamma _0u=\varphi \text{ on
}X',\tag2.9
$$
is solved by a a row matrix
$$
\pmatrix P_1-\lambda \\ \gamma _0
\endpmatrix^{-1}=
\pmatrix R_\lambda   & K_\lambda
\endpmatrix = \pmatrix Q_{\lambda ,+}+G_\lambda  & K_\lambda
\endpmatrix ,\text{ with }G_\lambda  =- K_\lambda  \gamma
_0Q_{\lambda ,+}.\tag2.10
$$
This is seen as follows: We denote by $K_\lambda $ the
Poisson operator solving the semi-ho\-mo\-ge\-ne\-ous problem with
$f=0$, i.e.,  $K_\lambda \varphi$ is the solution $u$ of the equations
$$(P_1-\lambda )u=0\text{ on }X,\quad \gamma _0u=\varphi \text{ on }X'.
$$
The resolvent $R_\lambda $ of $P_{1,\operatorname{D}}$ should solve the
other semi-homogeneous problem, with $\varphi =0$.
With $Q_{\lambda ,+}$
denoting the truncation to $X$ of the resolvent $Q_\lambda $ of $P_1$
on $\widetilde X$, one can easily check that
$R_\lambda =Q_{\lambda ,+}- K_\lambda  \gamma_0Q_{\lambda ,+}$ solves
the problem; it has the singular Green operator part
$G_\lambda  =- K_\lambda  \gamma_0Q_{\lambda ,+}$.

We shall also study powers of the resolvent and write them in the form
$$
R^N_\lambda =(Q_{\lambda ,+}+G_\lambda)^{N}=(Q^N_\lambda
)_++G^{(N)}_\lambda
.\tag 2.11
$$
Note that
$$R_\lambda ^{N}=\tfrac 1{(N-1)!}\partial _\lambda
^{N-1}R_\lambda ,\quad Q_{\lambda ,+}^{N}=\tfrac 1{(N-1)!}\partial
_\lambda
^{N-1}Q^N_{\lambda ,+},\quad G^{(N)}_\lambda =\tfrac
1{(N-1)!}\partial _\lambda ^{N-1}G_\lambda ;
\tag 2.12$$
this allows us to replace calculations for powers of $R_\lambda $ by
calculations for $\lambda $-derivatives.

We shall now start with the proof of Theorem 2.1. The case
where $A=G$ is a singular Green operator of class zero and
the pseudodifferential part vanishes can be attacked more
directly and will be studied first.

\head 3. Traces on the algebra of singular Green operators
\endhead

Let $G$ be a singular Green operator in $E$ of order $\nu \in\Bbb R$
and class 0,
and consider $GR^N_\lambda = GQ^N_{\lambda ,+}+GG^{(N)}_\lambda$.
In \cite{GSc} we introduced the auxiliary variable
$\mu =(-\lambda )^{\frac12}$ and phrased the results in terms
of $\mu$ instead of $\lambda $. In the present paper we shall
keep $\lambda $ as an index, but will often use $\mu $ in symbol
computations.
Some formulas from \cite{GSc} and some immediate consequences are
collected in Appendices A and B, to which we refer in the following.

As mentioned already,
the proof in \cite{GSc} of the expansion (1.10) for $A=G$ is
straightforwardly
modified to allow $\nu \notin\Bbb Z$.
For an analysis of the coefficients,
 we apply a partition of unity $1=\sum_{i=1}^{i_0}\theta _i$
subordinate to a cover of $X$ by coordinate patches $U _j$,
$j=1,\dots,j_0$, with trivializations $\psi _j\:E|_{U_j}\to V_j\times
\Bbb C^{\operatorname{dim}E}$,  $V_j\subset\subset \crnp$, $U_j\cap
\partial X$ mapped into $\partial \crnp$, such that any
two of the functions $\theta _{i_1}$ and $\theta _{i_2}$ are
supported in one of the coordinate patches $U _{j(i_1,i_2)}$.
Replacing $G$ by $\sum_{i_1,i_2\le i_0}\theta _{i_1}G\theta _{i_2}$,
we find from  $GQ^N_{\lambda ,+}$ and $GG^{(N)}_\lambda
$ a system of terms  $$
G=\sum_{i_1,i_2\le i_0}(\theta _{i_1}G\theta _{i_2}Q^N_{\lambda ,+}+
\theta _{i_1}G\theta
_{i_2}G^{(N)}_\lambda ),\tag3.1
$$
where $\theta _{i_1}$ and $\theta _{i_2}$ are supported either
in an interior coordinate patch (with closure in the interior) or a
patch meeting the boundary, and
we can study each term in the local coordinates.
Since $\theta _{i_2}G^{(N)}_\lambda $ is strongly polyhomogeneous of
order $-\infty $ when $\theta _{i_2}$ is supported in the interior,
its kernel is smooth and is $O(\lambda ^{-M})$ for
any $M$.
Hence the terms
$\theta _{i_1}G\theta_{i_2}G^{(N)}_\lambda $
that are supported in the interior have traces
that are $O(\lambda ^{-M})$ for any $M$ and can be disregarded in
the following.
For  the terms
$\theta _{i_1}G\theta _{i_2}Q^N_{\lambda ,+}$
supported in the
interior, $\theta _{i_1}G\theta _{i_2}$ is of order $-\infty $; they
will be covered by Lemma 3.1 
below, which gives:
$$
\Tr(\theta _{i_1}G\theta
_{i_2}Q^N_{\lambda ,+})=
\Tr(\theta _{i_1}G\theta
_{i_2})(-\lambda )^{-N}+O(\lambda ^{-N-1/2}).\tag3.2
$$
For the terms supported in patches meeting the boundary we can use
the results of \cite{GSc} in an accurate way.

Let us first recall some elements of the basic symbol calculus we
use, namely
the calculus of parameter-dependent $\psi $do's introduced in \cite{GS1}.
It will be used both in the \linebreak$(n-1)$-dimensional setting,
relevant for operators on $X'$, and in the $n$-dimensional setting,
relevant for operators on $\widetilde X$.

For $m\in\Bbb R$, the usual pseudodifferential symbol space $S^m(\Bbb
R^{n'}\stimes\Bbb R^{n'})$
consists of the functions $s(x,\xi )\in
C^\infty (\Bbb R^{n'}\stimes\Bbb R^{n'})$ satisfying estimates$$
|D_{x}^\beta D_{\xi }^\alpha s(x,\xi )|\le C_{\alpha ,\beta
}\ang{\xi }^{m-|\alpha |},
$$
for all $\alpha ,\beta \in \Bbb N^{n'}$. (We denote $(1+|\xi
|^2)^{\frac12}=\ang\xi $, $\{0,1,2,\dots\}=\Bbb N$.)
The parameter-dependent version we shall use here is $S^{m,d,s}(\Bbb
R^{n'}\stimes\Bbb R^{n'},\Gamma  )$, for $m\in\Bbb
R$, $d$ and $s\in\Bbb Z$, and with $\Gamma $ denoting a sector in
$\Bbb C\setminus \{0\}$; it was introduced in \cite{G3}, where a
detailed account can be found.
 Here $S^{m,d,s}(\Bbb
R^{n'}\stimes\Bbb R^{n'},\rp )$ is defined as the space of $C^\infty
$ functions $f(x,\xi ,\mu )$ on $\Bbb R^{n'}\stimes\Bbb
R^{n'}\stimes\rp $ such that, with $z=1/\mu $,$$
\multline
\partial _z^j[z^{d}|(\xi ,1/z)|^{-s}f(x,\xi ,1/z)]\in
S^{m+j}(\Bbb R^{n'}\stimes \Bbb R^{n'}),\\
\text{with symbol estimates uniform in $z$ for }z\le 1.
\endmultline$$
For more general $\Gamma $, the estimates have to hold on each ray in
$\Gamma $, uniformly in closed subsectors of $\Gamma $. When the
symbols moreover are holomorphic in
$\mu \in\Gamma ^\circ$ (just for $|(\mu ,\xi
')|\ge \varepsilon $ with some $\varepsilon $ depending on the closed
subsector), we speak of holomorphic symbols. The indication $\Bbb
R^{n'}\stimes \Bbb R^{n'}$ (and $\Gamma $) will often
be omitted. We say that $f(x,\xi
,\mu )$ is (weakly) polyhomogeneous in $S^{m,d,s}(\Gamma  )$,
when there is a
sequence of symbols $f_{j}$ in $S^{m-j,d,s}(\Gamma )$,
homogeneous in $(\xi ,\mu )$ of degree $m-j+d+s$ for $|\xi |\ge 1$, such
that
$f-\sum_{j<J}f_j$ is in $S^{m-J,d,s}(\Gamma )$ for all
$J\in\Bbb N$ (strongly polyhomogeneous if the homogeneity holds for
$|(\xi ,\mu )|\ge 1$, with appropriate remainder estimates).

These spaces $S^{m,d,s}$ are a generalization of the spaces $S^{m,d}$
introduced in \cite{GS1} (in fact, $S^{m,d}(\Gamma )$ equals the
space of holomorphic symbols in $S^{m,d,0}(\Gamma )$); they satisfy $$
\aligned
S^{m,d,s}&\subset S^{m+s,d,0}\cap
S^{m,d+s,0}\text{ if }s\le 0,\\
S^{m,d,s}&\subset S^{m+s,d,0}+
S^{m,d+s,0}\text{ if }s\ge 0,
\endaligned\tag3.3
$$
by \cite{GS1, Lemma 1.13}. The notation with the third upper
index spares us from keeping track of such
intersections and sums in general calculations. (It was used
systematically
in \cite{G3}, \cite{G5}, \cite{G6}.)

When $f\in S^{m,d,0}$, there is an expansion in decreasing powers of
$\mu $ starting with $\mu ^d$:
$$
f(x,\xi ,\mu )=\sum_{0\le k<N}\mu ^{d-k}f_{(k)}(x,\xi)+O(\mu
^{d-N}\ang\xi ^{m+N}),\quad \text{any $N$};\tag3.4
$$
in fact the $N$'th remainder is in $S^{m+N,d-N,0}$. Such expansions
play an important role in the proof of the following fact:
When
$F=\operatorname{OP}(f(x,\xi ,\mu ))$ has a polyhomogeneous symbol
in $S^{m,d,0}$ such that all terms
are integrable with respect to $\xi $, then by \cite{GS1, Th\. 2.1},
the kernel of $F$ has an asymptotic expansion on the diagonal:
$$
K(F,x,x,\mu )\sim \sum_{j\in\Bbb N}c_j(x)\mu ^{m+d+n'-j}+\sum_{k\in\Bbb
N}(c'_k(x)\log\mu +c''_k(x))\mu ^{d-k};\tag 3.5
$$
here the coefficients $c_j(x)$ and $c'_k(x)$ with $k=-m+j-n'$ are
local, determined from the strictly homogeneous version of the $j$'th
homogeneous term in the symbol of $F$, whereas the
$c''_k(x)$ depend on the full operator (are global). If
$m\notin\Bbb Z$, all $c'_k$ are zero.
When $F$ is a $\psi $do family on $\widetilde X$ (or on $X'$), the trace
expansion for
$F$ is found by establishing expansions (3.5) for components of $F$
in local coordinate
systems, carrying them back to $\widetilde X$ and integrating the
fibrewise traces  over
$\widetilde X$; this
leads to a similar
expansion of $\Tr F$:
$$
\Tr(F )\sim \sum_{j\in\Bbb N}c_j\mu ^{m+d+n'-j}+\sum_{k\in\Bbb
N}(c'_k\log\mu +c''_k)\mu ^{d-k}.\tag 3.6
$$
Of course, such a calculation depends on the choice of partitions of
unity and local coordinates, but it is useful for qualitative information.
The final result (3.6) is independent
of the choices, in the sense that the collected coefficient of each
power or log-power is so (for $k=j-m-n'\in\Bbb N$, $c_j+c''_k$ is
independent of the choices).

This is the basis for the trace expansions we shall show. For purely
$\psi $do terms on $X$, the expansions will be established relative to
$\widetilde X$, but the integration of the kernels will be restricted to
$X$. For all the other types of operators, the main
idea is to reduce the considerations to $\psi $do's over the
boundary $X'$. In all cases the point is to show that the resulting
symbols belong to suitable
$S^{m,d,0}$-spaces such that (3.5) can be established.

Since we are particularly interested in the nonlocal
coefficient $c_0''$, we have to be very careful with remainders and
smoothing terms.
Recall from \cite{GS1, Prop\. 1.21} that when an operator $F$ with
symbol in $S^{m,d,0}$
is of order $-\infty $ (i.e., lies in $S^{-\infty ,d,0}$), then it
has a $C^\infty $ kernel
$K(F,x,y,\mu )$ with an expansion in $C^\infty $ terms$$
K(F,x,y,\mu )=\sum_{0\le l<L}K_l(x,y)\mu ^{d-l}+K'_L(x,y,\mu ),\tag3.7
$$
for all $L$,
with $K'_L$ and its $(x,y)$-derivatives being $O(\ang \mu
^{d-L})$. 
Moreover, since $\partial _z^j(z^d F)=(-\mu ^2\partial _\mu )^j(\mu
^{-d}F)$ has symbol in $S^{-\infty ,0,0}$ for all $j$, there are expansions
$$
K((-\mu ^2\partial _\mu )^j(\mu
^{-d}F),x,y,\mu )=\sum_{0\le l<L}K^{(j)}_l(x,y)\mu ^{-l}+K^{(j)\prime}_L(x,y,\mu ),\tag3.8
$$
for all $L$,
with $K^{(j)\prime}_L$ and its $(x,y)$-derivatives being $O(\ang \mu
^{-L})$.
Conversely, these kernel properties characterize the operators $F$
with symbol in $S^{-\infty ,d,0}$ (as is seen straightforwardly from
the definition). 

It follows from (3.3) that $S^{-\infty ,d,s}=S^{-\infty ,d+s,0}$. 
In the proof that a polyhomogeneous symbol lies in a symbol space
$S^{m,d,s}$ one
shows the symbol properties for the homogeneous terms and treats the
remainders by establishing
the kernel expansionss (3.8) (with
$d$ replaced by $d+s$) with smoother coefficients and a higher number
of terms, the higher the index of the remainder is taken.

The sector $\Gamma $ used in the present work is$$
\Gamma =\{\,\mu \in\Bbb C\setminus \{0\}\mid |\arg \mu |\le\tfrac \pi
4+\tfrac
\varepsilon 2
\,\},\tag3.9
$$
for some $\varepsilon>0$,
since $\lambda =-\mu ^2$ runs in $\Lambda $ defined in (2.3).
{\it All symbols to be considered in this paper will be holomorphic
on $\Gamma ^\circ$}; this fact will not be mentioned explicitly each time.

One of the technical points in \cite{GSc} was that in compositions
of two operators, the right-hand factor is taken in
$y$-form and the resulting symbol is considered in $(x,y)$-form so
that one does not need to consider the asymptotic expansion in the
usual symbol composition formula (cf.\ e.g.\  (3.31) below). However,
the passage from one of these
forms to another can induce low-order errors. These
have to be dealt with
carefully in the present situation, where nonlocal contributions are
in focus. Also for this reason, the following analysis of the effect
of smoothing operators is important.

Smoothing parts of $Q_\lambda $ and $G_\lambda $ are easy to deal with,
since they have smooth kernels that are $O(\lambda ^{-M})$ for any
$M$.
For compositions containing smoothing parts of $G$, we recall that an
s.g.o\. of order
$-\infty $ is simply an operator with $C^\infty $ kernel on $X$.  So is
$P_+$ when $P$ is of order $-\infty $. Such
operators  enter as follows (for $G^{\pm}$, cf\. (A.14)):

\proclaim{Lemma 3.1} Let $G$ be an operator with $C^\infty $ kernel,
compactly supported in $\crnp\stimes\crnp$. Then
the traces of $GQ^N_{\lambda ,+}$, $GG^{(N)}_\lambda $ and
$GG^\pm(Q^N_\lambda )$ have expansions
$\sum_{j\ge 0}c''_j\mu ^{-2N-j}$, where
$$
\alignedat2
\text{\rm (i)}&&\quad
c''_0&= \Tr G \text{ in the expansion of }\Tr(GQ^N_{\lambda ,+}),\\
\text{\rm (ii)}&&\quad c''_0&= 0 \text{ in the expansions of
$GG^{(N)}_\lambda $ and
$GG^\pm(Q^N_\lambda )$.}
\endaligned \tag3.10
$$
\endproclaim

\demo{Proof} First consider the cases in (ii); here
we use that the symbol of the $\lambda
$-dependent factor has the structure described in Lemma A.4
with $\mu^2=-\lambda$.
Consider a term $G_{J,j,j'}=\operatorname{OPG}(g_{J,j,j'}(x',\xi ',\xi
_n,\eta _n,\mu ))$ with $g_{J,j,j'}=s_{J,j,j'}(\kappa +i\xi
_n)^{-j}(\kappa '-i\eta _n)^{-j'}$ ($\kappa $ and $\kappa '$
taking values $\kappa ^+$ or $\kappa ^-$), $s_{J,j,j'}\in
S^{0,0,j+j'-1-2N-J}$; its kernel is
$$\aligned
K(G_{J,j,j'},x,y,\mu )
&=
\int e^{i (x'-y')\cdot \xi '+ix_n\xi
_n-iy_n\eta _n}g_{J,j,j'}(x',\xi ',\xi _n,\eta _n,\mu )\,\d\xi \d\eta
_n\\
&= c_jc_{j'}
\int e^{i (x'-y')\cdot \xi '}s_{J,j,j'}(x',\xi ',\mu
)x_n^{j-1}y_n^{j'-1}e^{-\kappa x_n-\kappa 'y_n}\,\d\xi ',
\endaligned
$$
since $\Cal F^{-1}_{\xi _n\to x_n}(\kappa +i\xi
_n)^{-j}=c_jx_n^{j-1}e^{-\kappa x_n}$.
Now $$\aligned
K(GG_{J,j,j'},x,y,\mu )&=\int_{\rnp}K(G,x,z)K(G_{J,j,j'}, z,y,\mu )\,dz
,\\
\Tr (GG_{J,j,j'})&=\int_{\rnp}\tr\int_{\rnp}K(G,x,z)K(G_{J,j,j'}, z,x,\mu
)\,dzdx .\endaligned
$$
To evaluate the latter, we insert a Taylor expansion of $K(G,x,z)$ w.r.t.\
$(x_n,z_n)$:$$
K(G,x,z)=\sum_{m,m'<M}b_{m,m'}(x',z')x_n^mz_n^{m'} +x_n^Mz_n^Mb'_M(x,z),
$$
and use that $$
\aligned
\int_0^\infty x_n^{m+j-1}e^{-\kappa x_n}\,dx_n&=c'_{m+j}\kappa
^{-m-j},\\
|\int_0^\infty f(x_n)x_n^{M+j-1}e^{-\kappa x_n}\,dx_n|&\le
\sup |f(x_n)e^{-\kappa x_n/2}\int_0^\infty x_n^{M+j-1}e^{-\kappa
x_n/2}\,dx_n|\\
&\le
C\kappa
^{-M-j}.
\endaligned
$$
This gives when
integrations in $(x_n,z_n)$ are performed:
$$\multline
\Tr (GG_{J,j,j'})=\sum_{m,m'<M}c_{jj'mm'}\\
\int\tr
b_{m,m'}(x',z')\kappa ^{-m-j'}(\kappa ')^{-m'-j}e^{i(z'-x')\cdot \xi
'}s_{J,j,j'}(
z',\xi ',\mu )\,\d \xi 'dz'dx'\\
\quad+\int O(\ang{(\xi ',\mu )} ^{-2M-1-2N-J})\,\d\xi '.
\endmultline$$
Each term in the sum over $m,m'$ is the trace of a $\psi $do
$S_{J,j,j',m,m'}$ on $\Bbb
R^{n-1}$ composed of the operator with kernel $c_{jj'mm'}b_{m,m'}(x',y')$
(smooth compactly supported) and
the operator with symbol $$\kappa (x',\xi ',\mu )^{-m-j'}\kappa
'(x',\xi ',\mu )^{-m'-j}s_{J,j,j'}(
x',\xi ',\mu )\in S^{0,0,-m-m'-1-2N-J};
$$
so $S_{J,j,j',m,m'}$ has symbol in $S^{-\infty
,0,-m-m'-1-2N-J}=S^{-\infty ,-m-m'-1-2N-J,0}$
and consequently a trace expansion as in (3.7) with $d\le-2N-1$ in
all cases. The remainder is $O(\mu ^{-2M-1-2N-J+n-1})$. Since $M$ can
be taken arbitrarily large, this shows (ii).

For (i), let $\widetilde G$ be an operator on $\Bbb R^n$ with
compactly supported $C^\infty $ kernel such that $G=r^+\widetilde G
e^+$ (its kernel can be constructed by extending the kernel of $G$
smoothly across $x_n=0$ and $y_n=0$). Then $$
GQ^N_{\lambda ,+}=\widetilde G_+Q^N_{\lambda ,+}=(\widetilde
GQ^N_{\lambda })_+-G^+(\widetilde G)G^-(Q^N_{\lambda }).
$$
Here $G^+(\widetilde G)G^-(Q^N_{\lambda })$ is of the type we
considered above, having a trace expansion with $d\le -2N-1$.
For $(\widetilde
GQ^N_{\lambda })_+$ we use the result known from
\cite{GS1} (cf.\ also \cite{G4, Th\. 1.3})
that
the kernel of $\widetilde GQ^N_{\lambda }$ has an expansion on the
diagonal $$
K(\widetilde GQ^N_{\lambda },x,x,\mu )\sim \sum_{j\ge 0}c''_j(x)\mu
^{-2N-j},\text{ with }c''_0(x)=K(\widetilde G,x,x ).\tag 3.11
$$
 The expansion of
$\Tr(\widetilde GQ^N_\lambda )_+$ is found by integrating the matrix
trace of (3.11)
over $\rnp$, where $K(\widetilde G,x,x )=K( G,x,x )$, so (i)
follows, with $c''_0=\int_{\rnp}\tr K( G,x,x )\, dx=\Tr G$.
\qed
\enddemo

As already mentioned, this lemma takes care of the interior parts of
$GQ^N_{\lambda ,+}$, cf\. (3.1), (3.2).

For the parts $\theta
_{i_1}G\theta
_{i_2}G^{(N)}_\lambda $ at the boundary, we can use Proposition B.3
directly;
it shows that they contribute only locally to the coefficient of
$\lambda^{-N}$.

It remains to consider terms $\theta _{i_1}G\theta
_{i_2}Q^N_{\lambda ,+}$ in patches intersecting the boundary.
Here we can draw on the fine analysis of such terms made in
\cite{GSc} to find
 the
contribution to $c'_0$ in (1.12); it can be used also to isolate
the contribution to $c''_0$ modulo local terms. This will give not
only a proof of (2.7) but also a ``value'' (modulo local terms)
of the constant.

To simplify the presentation, we can unify the treatment of the
coordinate patches by spreading out the images $V_{j(i_1,i_2)}$ of
the sets $U_{j(i_1,i_2)}$ by linear translations in the $x'$-variable,
to obtain sets $V'_{i_1,i_2}$ with positive distance (in the
$x'$-direction) from one another. Then the sum of
the localized operators $\theta _{i_1}G\theta _{i_2}$ acts in
$\crnp\times \Bbb C^{\operatorname{dim}E}$; we shall denote it by $G$
again or, if a distinction from the original $G$ is needed, by $\underline
G$. $R_{\lambda }$ is likewise considered in these coordinates and
may be denoted $\underline R_\lambda $ if needed for precision.

Recall
that when $G$ is defined on $\crnp$ from a symbol $g(x',\xi ',\xi
_n,\eta _n)$ with Laguerre expansion
(B.15), then the operator $\tr_nG$ is defined as the $\psi $do on $\Bbb
R^{n-1}$ with symbol
$$\tr_ng(x',\xi ')=\int g(x',\xi ',\xi _n,\xi
_n)\,\d\xi _n=\sum_{j\in\Bbb N}d_{jj}(x',\xi ').$$
Recall also (cf\. e.g\.
\cite{G1}) that a singular Green operator on $\overline {\Bbb
R}^n_+$ with compact $x'$-support is trace-class, when its
order is $<-n+1$. Then $\tr_nG$ is trace-class on $\Bbb R^{n-1}$ and
$\Tr _{\Bbb R^n} G=\Tr_{\Bbb R^{n-1}}(\tr_n G)$.

We shall now apply the result of  Proposition B.5; that
$\tr_n(GQ^N_{\lambda ,+})=\widetilde S_0+\widetilde S_1$, with
$\widetilde
S_0=\op'(\tr_ng(x',\xi
')\alpha ^{(N)}(x',\xi ',\mu ))$ (cf\. (A.12)) and $\widetilde S_1$
having symbol in $S^{\nu +1,-2N-1,0}\cap S^{\nu -2N,0,0}$. This
result was used in  \cite{GSc} to pinpoint the first
logarithmic coefficient $\tilde c_0'$ in (1.10), using that only
$\widetilde S_0$ contributes to it. But the fact that $\widetilde
S_1$  has $d$-index $-1-2N$ also implies that
it contributes only locally  to $\tilde c_0''$, so nonlocal contributions
to $\tilde c''_0$ come entirely from $\widetilde S_0$.

In the localized situation, denote $\tr_nG=\ttilde G$ and denote its
symbol $\tr_ng=\ttilde g$, expanded in
homogeneous terms $\ttilde g\sim\sum_{j\ge 0}\ttilde g_{\nu -j}$.

We shall use the notion of a
regularized integral (or finite part integral) $\tslint f(\xi
)\,\d\xi $ as in Lesch
\cite{L}, \cite{G4}, for polyhomogeneous  functions
$f(\xi )$ of $\xi \in\Bbb R^{n'}$
(it will be used with $n'=n-1$ or $n$):
$\tslint f( \xi )\,\d\xi $ is the constant term in the expansion
of $\int_{|\xi|\le\mu }f(\xi) \d\xi$ 
into powers of $\mu$ and $\log \mu$.
In more detail:

When $f( \xi )$ is  integrable in $\xi $,
$\tslint f( \xi )\,\d\xi $ is the usual integral $\int_{\Bbb
R^{n'}}f( \xi )\,\d\xi $. When $f_{\nu-j}( \xi )$ is
homogeneous of order
$\nu -j$ in $\xi $ for $|\xi |\ge 1$, then
$$
\int_{|\xi |\le \mu }f_{\nu -j}\,\d\xi =\int_{|\xi |\le 1 }
f_{\nu -j}\,\d\xi +\int_{1\le |\xi |\le \mu  }f_{\nu -j}\,\d\xi ,\tag3.12
$$
where a calculation in polar coordinates gives
$$
\aligned
\int_{1\le |\xi |\le \mu }&f_{\nu -j}( \xi)\,\d\xi
=\int_1^\mu r^{\nu -j +n'-1}\,dr
\int_{|\xi |=1}f_{\nu -j} ( \xi )\,\d S(\xi )\\
&=\cases 
c_{j}(\frac{\mu ^{\nu -j +n'}}{\nu -j+n'}-\frac1{\nu
-j+n'})&\text{ if }\nu -j+n'\ne 0,\\
c_{j}\log \mu &\text{ if }\nu -j+n'= 0;\endcases\\
c_{j}&= \int_{|\xi |=1}f_{\nu -j}( \xi
)\,\d S(\xi ).\endaligned\tag3.13
$$
Here,
$$
\slint f_{\nu -j}( \xi )\,\d\xi =\int_{|\xi |\le 1}f_{\nu -j}\,\d\xi
+\cases -\frac1{\nu
-j+n'}c_j&\text{ if }\nu -j+n'\ne 0,\\
0 &\text{ if }\nu -j+n'= 0;\endcases\tag3.14
$$
this is consistent with the integrable case.
(Such calculations were
basic in the proof of \cite{GS1, Th\. 2.1}.)
More generally, when
$f(x, \xi )$
is a classical symbol of
order $\nu $, $f(x, \xi )\sim \sum_{j\in\Bbb N}f_{\nu -j}(x,\xi
)$, then there is an asymptotic expansion for $\mu \to\infty $,
$$
\int _{|\xi |\le \mu }f(x,\xi )\,\d\xi \sim \sum_{j\in\Bbb N, j\ne
\nu +n'}f _j(x)\mu ^{n'+\nu -j}+f '_0(x)\log \mu +f ''_0(x),
\tag3.15$$
and we set $$
\slint f(x,\xi )\,\d\xi =f ''_0(x).\tag3.16
$$
(This is related to Hadamard's
definition of the finite part --- {\it partie finie} --- of certain
integrals \cite{H, p\. 184ff.}.)
In view of (3.12)--(3.14), we have the precise formula:
$$
\multline
\slint f(x,\xi )\,\d\xi
=\sum_{0\le j\le \nu +n'}\bigl(\int_{|\xi|\le 1}
f_{\nu -j}(x,\xi )\,\d\xi -
\tfrac{1-\delta _{\nu +n',j}}{\nu +n'-j}
\int_{|\xi |=1}f_{\nu -j}\,\d S(\xi )\bigr)\\
+\int_{\Bbb R^{n'}}
(f-\sum_{j \le \nu +n'}f_{\nu -j})\,\d\xi ,
\endmultline\tag3.17
$$
where $\delta _{r,s}$ is the Kronecker delta. (One can replace the
sum over $j\le \nu +n'$ by the sum over $j\le J$ for any $J\ge \nu +n'$.)

\example{Remark 3.2}
Let $F=\op (f(x,\xi))$ be the $\psi$do associated with the above
symbol. It was shown by Lesch \cite{L, Section 5} that the
density $\omega (F)=\tslint f(x,\xi)\,\d\xi dx$ is invariant under
the change
of the symbol induced by diffeomorphisms of open sets, when
$\nu\not\in\Bbb Z$ or  $\nu<-n'$. Moreover, it was observed in \cite{G4}
that
the proof of \cite{L} extends to the cases where $\nu \in \Bbb Z$ and
$f$ has a parity that fits with the dimension $n'$, cf\.
(1.4)--(1.5)ff. So in all these cases, $\omega (F)$ depends only on
$F$, not on the representation of its symbol in
local coordinates. As a consequence, we can consider this density also
for $\psi$do's having the mentioned properties on a manifold $M$ of
dimension $n'$.
Note that when $M$ is compact, $\int_M\tslint \tr
f(x,\xi)\,\d\xi dx=\Tr(F)$ if $\nu<-n'$.

\endexample

The notation is used in the following with $n'=n-1$.
It was shown in \cite{G4, Th\. 1.3} (by working out
\cite{GS1, Th\. 2.1} explicitly in this case) that when
$f(x',\xi ')$ is the symbol of a classical $\psi $do $F$ on $\Bbb
R^{n-1}$, and $S=\operatorname{OP}'(s(x',\xi '))$ is an auxiliary
second-order uniformly elliptic operator with no principal
eigenvalues on $\Bbb
R_-$, then the kernels of $F(S+\mu
^2)^{-N}$
and $\operatorname{OP}'(f(x',\xi ')(s(x',\xi
')+\mu ^2)^{-N})$ with $N>(\nu +n-1)/2$ have diagonal expansions
of the form$$
K(x',x',\mu )\sim \sum_{j\ge 0}b_j(x')\mu ^{n-1+\nu
-j-2N}+\sum_{k\ge 0}(b'_k(x')\log \mu +b''_k(x'))\mu
^{-k-2N}, \tag3.18
$$ where the $b'_k$ are 0 if $\nu \notin\Bbb Z$, and
$$
\aligned
b_{\nu +n-1}(x')+b''_0(x')&=\slint f(x',\xi ')\,\d\xi '+\text{ local
terms, if $\nu $ is integer $\ge 1-n$},\\
b''_0(x')&=\slint f(x',\xi ')\,\d\xi ' \text{ if $\nu <1-n$ or }
\nu \notin \Bbb Z.\endaligned
\tag3.19
$$
Moreover,
if $\nu \in\Bbb Z$ and $f$ has a parity that
fits with the dimension $n-1$, then $b_{\nu +n-1}(x')=0$ and the second
formula in (3.19) holds.
(Furthermore, the sum over $k$ in (3.18) for the operators $F(S+\mu
^2)^{-N}$
and $\operatorname{OP}'(f(s+\mu ^2)^{-N})$ skips the
terms where $k/2\notin \Bbb Z$, but we need to refer to the general
expansion below. In parity cases, the $b_j$ are zero for $\nu +n-1-j$
even, and the $b'_k$ are zero.)  
What we shall show now is that $\alpha ^{(N)}$ plays to a large extent the
same role as $(s+\mu ^2)^{-N}$.
We shall deal with the special cases where $\nu <1-n$ or
$\nu \notin \Bbb Z$ here, whereas parity cases will be discussed
later, around Theorem 3.15.

\proclaim{Theorem 3.3} One has for $\widetilde
S_0=\operatorname{OP}'(\ttilde g\alpha ^{(N)})$,
when $N>(n-1+\nu
)/2$: The kernel of
$\widetilde S_0$ has an
expansion on the diagonal:
$$
K(\widetilde S_0,x',x',\mu )
\sim \sum_{j\ge 0} a_{j}(x')\mu ^{{n-1+\nu - j}-2N}+\sum_{k\ge
0}(a'_k(x')\log \mu + a''_k(x'))\mu
^{- k-2N},\tag 3.20
$$
where the terms $a_j(x')$, and $a'_k(x')$ for $k=\nu +n-1-j$, depend on
the
first $j+1$ homogeneous terms in the symbols $\ttilde g(x',\xi ')$ and
$\alpha ^{(N)}(x',\xi ',\mu )$; the $a''_k(x')$ are global.

If $\nu <-n+1$ or $\nu \in \Bbb R\setminus \Bbb Z$, then
$$
a''_0(x')=\slint
\ttilde g(x',\xi ')\,\d\xi ',
\tag3.21
$$
whereas if $\nu $ is an integer $\ge -n+1$,$$
a_{\nu +n-1}(x')+a''_0(x')=\slint
\ttilde g(x',\xi ')\,\d\xi '+\text{ local terms},
\tag3.22
$$
where the local terms depend on the first $\nu +n$ homogeneous terms
in the symbols of $\ttilde g$ and $\alpha ^{(N)}$. The $a'_k$ vanish
if $\nu \notin\Bbb Z$.

\endproclaim

\demo{Proof} The existence of the expansion is assured by \cite{GS1,
Th\. 2.1}, since $\widetilde S_0$ has a weakly
polyhomogeneous symbol in $S^{\nu ,-2N,0}\cap S^{\nu -2N,0,0}$, so
the main point is to prove the
formulas for  $a''_0$ resp.\ $a_{\nu +n-1}+a''_0$.
Recall that $$
K(\widetilde S_0,x',x',\mu )=\int_{\Bbb R^{n-1}}\ttilde g(x',\xi
')\alpha ^{(N)} (x',\xi ',\mu )\,\d\xi '.
$$
First consider
$\ttilde g'=\ttilde g-\sum_{j\le
\nu +n-1}\ttilde g_{\nu -j}$, it is of order $\nu -J<-n+1$ where
$J=\max\{n+[\nu ],0\}$, hence
integrable in $\xi '$. Since$$
\alpha ^{(N)}=\mu ^{-2N}+\alpha ^{(N)}_1(x',\xi ',\mu )
$$ with $\alpha ^{(N)}_1\in S^{0,-2N,0}\cap S^{1,-2N-1,0}$
(cf\. (A.13)ff.),
$\ttilde g'\alpha ^{(N)}_1$ is in $S^{\nu -J,-2N,0}\cap
S^{1+\nu -J,-2N-1,0}$,
and we can write
$$
K(\operatorname{OP}'(\ttilde g'\alpha ^{(N)} ),x',x',\mu )=\mu
^{-2N}\int_{\Bbb R^{n-1}}\ttilde
g'\,\d\xi '+\int_{\Bbb R^{n-1}}\ttilde
g'\alpha ^{(N)}_1\,\d\xi '.\tag3.23
$$
By \cite{GS1, Th\. 2.1}, the last integral has an expansion as in
(3.18) with $\nu $ replaced by $\nu -J$ and the sum over
$k$ starting with $k=1$. This expansion has no term with $\mu
^{-2N}$,
so the only contribution from
 $\ttilde g'$ to the coefficient of $\mu
^{-2N}$ is $\int \ttilde g' \,\d\xi =\tslint \ttilde g' \,\d\xi$.
This ends the proof in the case $\nu <-n+1$, where $\ttilde g=\ttilde
g'$.

For the homogeneous terms $\ttilde g_{\nu -j}$, we use an analysis
as in \cite{GS1, Th\. 2.1}:
Let $\mu \in\Bbb R_+$ and write the contribution
from $\ttilde g_{\nu -j}$ as $$
\int_{|\xi '|\ge \mu }\ttilde g_{\nu -j}(x',\xi
')\alpha ^{(N)} (x',\xi ',\mu )\,\d\xi '+\int_{|\xi '|\le \mu }\ttilde
g_{\nu -j}(x',\xi
')\alpha ^{(N)} (x',\xi ',\mu )\,\d\xi '.
$$
The first integral gives a local term $\tilde a_j(x')\mu ^{n-1+\nu
-j-2N}$, since
the integrand is homogeneous in $(\xi ',\mu )$ of degree $\nu -j-2N$.
Note that when $\nu \in\Bbb R\setminus \Bbb Z$, the power cannot be
${-2N}$.
The second integral may be written$$
\int_{|\xi '|\le \mu }\ttilde g_{\nu -j}\alpha ^{(N)} \,\d\xi '=\mu
^{-2N}\int_{|\xi '|\le \mu }\ttilde g_{\nu -j}(x',\xi
')\,\d\xi '+\int_{|\xi '|\le \mu }\ttilde g_{\nu -j}\alpha^{(N)}_1
\,\d\xi '.\tag3.24
$$
Here $\int_{|\xi '|\le \mu }\ttilde g_{\nu -j}\,\d\xi '$ is
calculated as in (3.12)--(3.14), contributing
$\tslint \ttilde g_{\nu -j}\,\d\xi '$ to the coefficient of $\mu
^{-2N}$.
For the last integral in (3.24) we observe that since $\ttilde
g'\alpha ^{(N)}_1$ is in \linebreak$S^{\nu -J,-2N,0}\cap
S^{1+\nu -J,-2N-1,0}$, with $d$-index $\le -2N-1$, the consideration
of this integral in the proof of \cite{GS1, Th\. 2.1} shows that it
contributes only locally to the coefficient of $\mu ^{-2N}$, and not at
all if $\nu \notin\Bbb Z$.
\qed
\enddemo

\example{Remark 3.4}
Similar statements hold when $\ttilde g$ is given in $(x',y')$-form;
then $\ttilde g(x',\xi ')$ is replaced by $\ttilde g(x',x',\xi ')$ in
{\rm (3.21)} and {\rm (3.22)}. The extension to this case is carried
out as in \cite{G4, Remark 1.4}.

\endexample

In view of the information on general compositions recalled before
the theorem, we have in particular:

\proclaim{Corollary 3.5} Let $S=\operatorname{OP}'(s(x',\xi '))$ be a
second-order uniformly elliptic operator on $\Bbb R^{n-1}$ with no
principal
eigenvalues on $\Bbb
R_-$. Comparing
the expansion {\rm (3.18)} for the kernel of \linebreak
$\ttilde G(S+\mu ^2)^{-N}$
with the expansion {\rm (3.20)} for the kernel of
$\widetilde S_0$,
one has that
$$
\aligned
a''_0(x')&=b''_0(x')+\text{ local
terms, if $\nu $ is integer $\ge 1-n$},\\
a''_0(x')&=b''_0(x')\text{ if $\nu <1-n$ or }
\nu \notin \Bbb Z.\endaligned
\tag3.25
$$
\endproclaim

We then find:

\proclaim{Theorem 3.6} Let $G$ be a singular Green operator on $X$
of order $\nu \in\Bbb R$ and class $0$.

$1^\circ$ Let $g(x',\xi ,\eta _n)$ be the symbol of the operator in a
localization to $\crnp$ as described above. Then
$$
C_0(G,P_{1,\operatorname{D}})=
\int\slint
\tr (\tr_ng)(x',\xi ')\,\d\xi 'dx'+\text{ local
terms};\tag 3.26
$$
it holds with
vanishing local terms if $\nu <-n+1$ or $\nu \in\Bbb R\setminus \Bbb Z$.

$2^\circ$ Let $P_1$ and $P_2$ be
two choices of auxiliary strongly elliptic operators. Then
$$C_0(G,P_{1,\operatorname{D}})- C_0(G,P_{2,\operatorname{D}})
\tag3.27$$ is
locally determined;
it  vanishes if $\nu <-n+1$ or $\nu \in \Bbb
R\setminus\Bbb Z$.
\endproclaim

\demo{Proof}
Consider $G$ and $Q^N_{\lambda ,+}$ carried over to $\crnp$.
Here we have
$$
\Tr(GQ^N_{\lambda ,+})=\Tr_{\Bbb
R^{n-1}}\tr_n(GQ^N_{\lambda ,+})=
\Tr_{\Bbb R^{n-1}}\widetilde S_{0}+\Tr_{\Bbb R^{n-1}}\widetilde S_{1},
\tag3.28$$
as in Proposition B.5. Since $\widetilde S_{1}$ has symbol in $S^{\nu
+1,-2N-1,0}\cap S^{\nu -2N,0,0}$, it has a diagonal kernel expansion as in
(3.18) with the sum over $k$ starting with $k=1$, so the coefficient
of $\mu ^{-2N}$ is locally determined (vanishing if $\nu \in\Bbb
R\setminus\Bbb Z $ or $\nu <1-n$). For $\widetilde S_{0}$, we have the
diagonal kernel expansion from Theorem 3.3. We integrate the matrix
traces  in $x'$, and we add on $\Tr(GG^{(N)}_\lambda )$, which is
described in Proposition B.3. This gives a trace expansion as in (1.10), with
coefficient of $(-\lambda )^{-N}=\mu ^{-2N}$ equal
to
$$
\int \slint \tr \tr_ng\,\d\xi 'dx'+\text{
local terms}
$$
(no local terms if $\nu \in\Bbb
R\setminus \Bbb Z $ or $\nu <1-n$); this shows $1^\circ$.

Now if $P_1$ is replaced by $P_2$ in these calculations, $1^\circ$
implies that
the contribution to the coefficient of $\mu ^{-2N}$ is modified only
in the local terms, vanishing if $\nu <-n+1$ or $\nu \in \Bbb
R\setminus\Bbb Z$; this shows $2^\circ$.
\qed
\enddemo

\example{Remark 3.7}
(a) Another choice of local coordinates and partition of unity will
give another decomposition
in (3.1), but this will
modify the term of a given order in the localized symbol only by
terms of the same and higher order.

(b) When $\nu \notin \Bbb Z$ or $\nu<-n+1$, or $\nu \in\Bbb Z$ and
$\tr_n g$ has a parity that fits with the dimension $n-1$, then
$\int \tslint \tr
\tr_ng\,\d\xi 'dx'$
will be
invariant under such choices, by Remark 3.2.

(c) In the localized situation, we find by integration from Corollary 3.5
that
$$
C_0( G, P_{1,\operatorname{D}})=C_0(\tr_n
 G, S)
+\text{ local terms,}\tag 3.29
$$
for any auxiliary operator  $S$ as described there. This reduces the
calculation of $C_0$ on $\rnp$ to a calculation of a $C_0$ on $\Bbb
R^{n-1}$, modulo local terms (vanishing if $\nu <1-n$
or $\nu \notin \Bbb Z$).
\endexample

We have hereby obtained (2.7) in the case $A=G$.
Now we turn to (2.8), the commutation property of
$C_0(G,P_{1,\operatorname{D}})$.
We shall give a proof that reduces it to the commutation property for
closed manifolds (cf\. e.g.\ \cite{G4}), using a variant of the
preceding analysis.

Here we want to consider $G$ given in a different form in the
localized situation,
namely in the form $$
G=\sum_{j,k\in\Bbb N}\Phi _jC_{jk}\Phi ^*
_k\tag 3.30$$
explained in Lemma B.6. When
$G=\sum_{j,k\in\Bbb N}\Phi _jC_{jk}\Phi ^*
_k$, then $$
\Tr_{\rnp}(GQ^N_{\lambda ,+})=\Tr_{\rnp}(\sum_{j,k\in\Bbb N}\Phi
_jC_{jk}\Phi ^*
_k Q^N_{\lambda ,+})=\sum_{j,k\in\Bbb N}\Tr_{\Bbb R^{n-1}}(C_{jk}\Phi ^*
_k Q^N_{\lambda ,+}\Phi _j),
$$
since the sum over $j,k$ is rapidly decreasing in the relevant symbol- and
operator norms. Here the $C_{jk}$ and $\Phi ^*
_k Q^N_{\lambda ,+}\Phi _j$ are $\psi $do's on $\Bbb R^{n-1}$; to use
their
properties, we have to investigate
$\Phi ^* _k Q^N_{\lambda ,+}\Phi _j$.

It was a typical feature of the calculations in
\cite{GSc} that we used the passage from $x'$-form to $y'$-form
and vice versa in the compositions in such a way that the
case of $\psi $do symbols over the boundary was reached before it was
necessary to
verify asymptotic composition formulas such as
$$
a(x',\xi ',t)\circ b(x',\xi ',s)\sim \sum_{\alpha \in \Bbb
N^{n-1}}\tfrac{1}{\alpha !}D _{\xi'} ^\alpha
a(x',\xi ',t)\partial^\alpha _{x '}b(x',\xi ',s),\tag3.31
$$ with parameters $t$ and $s$. In the $S^{m,d,s}$-spaces
over the boundary,
the composition formulas work as usual (\cite{GS1, Th\. 1.18}),
whereas the mixture of $\mu $-dependence and $\mu $-independence
makes it harder to see what goes on in compositions of the $\psi $dbo
symbols reaching
into  the interior of $X$. Here we rely heavily on the special
rational structure of the symbols entering in $Q_\lambda $. (The
compositions are
not covered by \cite{G3}.)

In compositions containing $\Phi ^*
_k Q^N_{\lambda ,+}\Phi _j$, we have to deal with terms of
the form$$
\operatorname{OPT}(\bar{\hat\varphi
}_k(\xi _n,\sigma (\xi ') )) \op (q_{-2-J}(x,\xi ,\mu
))_+\operatorname{OPK}(\hat\varphi _j(\xi _n,\sigma (\xi ') ))
$$
(and their $\lambda $-derivatives),
which do require asymptotic expansions. In fact, when
$q_{2-J}$ is independent of $x_n$, the two factors to the right
give a relatively simple Poisson symbol depending on $x'$, but then
in the composition with
the trace symbol to the left (depending on $\xi '$ through $\sigma (\xi
')$), there will be an infinite expansion as in (3.31), where
remainder estimates have to be shown. We shall treat this by going
back to the original
proof of the pseudodifferential composition rule by Taylor expansion,
keeping track of the
properties of the terms, in particular the remainder, by use of the
exact formulas. A composition to the left with
$C_{jk}=\op'(c_{jk}(x',\xi '))$ does not make
the expression harder to deal with, since $\op'(c_{jk}(x',\xi
'))\operatorname{OPT}(\bar{\hat\varphi
}_k(\xi _n,\sigma )) =\operatorname{OPT}(c_{j,k}(x',\xi ')\bar{\hat\varphi
}_k(\xi _n,\sigma ))$, so we include that in the following.

\proclaim{Proposition 3.8} Let $C_{jk}=\op'(c_{jk}(x',\xi '))$ be
of order $\nu $ and let
$$
S_{J,j,k}=\op'(c_{jk}(x',\xi '))\operatorname{OPT}(\bar{\hat\varphi
}_k(\xi _n,\sigma )) \op (q_{-2-J}(x,\xi ,\mu
))_+\operatorname{OPK}(\hat\varphi _j(\xi _n,\sigma )).\tag 3.32
$$
For each $J\in\Bbb N$,
$S_{J,j,k}$ is a $\psi $do with symbol
in $S^{\nu ,0,-2-J}$,
the $(N-1)$'st
$\lambda $-derivatives of the symbol lying in $S^{\nu ,0,-2N-J}$.

\endproclaim

\demo{Proof}
In view of the expansion in (A.5), $S_{J,j,k}=\sum_{J/2+1\le m\le
2J+1}S_{J,j,k,m}$, where
$$
\aligned
S_{J,j,k,m}&=C_{jk}\Phi _k^*\op \Bigl(\frac{r_{J,m}(x,\xi
)}{(p_{1,2}+\mu ^2)^{m}}\Bigr)_+\Phi _j\\
&=\operatorname{OPT}(c_{jk}(x',\xi ')\bar{\hat\varphi
}_k(\xi _n,\sigma )) \op \Bigl(\frac{r_{J,m}(x,\xi
)}{(p_{1,2}+\mu ^2)^{m}}\Bigr)_+\operatorname{OPK}(\hat\varphi _j(\xi
_n,\sigma )),
\endaligned\tag 3.33
$$
and we have to show the property for
each of the terms in the expansion.
The function $r_{J,m}$ is a polynomial
in $\xi $ of degree $2m-J-2\ge 0$. Let us drop the indexations on
$r_{J,m}$ and $p_{1,2}$, and simply write the fraction as $r(p+\mu
^2)^{-m}=q$.

We first consider the case where these symbols are independent of
$x_n$. Here one has immediately:
$$
\operatorname{OP}(q(x',\xi ,\mu
))_+\operatorname{OPK}(\hat\varphi _j(\xi _n,\sigma
))=\operatorname{OPK}(b(x',\xi ,\mu )),\quad b=h^+_{\xi _n}(q\hat\varphi
_j).
$$
To prepare for the composition of this Poisson operator with the trace
operator \linebreak$\operatorname{OPT}(c_{jk}\bar{\hat\varphi
}_k)$ to the left, we apply the usual procedure for changing an operator
from
$x'$-form to
$y'$-form:
$$
\multline
\operatorname{OPK}(b(x',\xi ,\mu ))v=\int e^{i(x'-y')\cdot \xi
'+ix_n\xi _n}b(x',\xi ,\mu )v(y')\,dy'\d \xi \\
=\int e^{i(x'-y')\cdot \xi
'+ix_n\xi _n}\sum_{|\alpha |<M}
\tfrac 1{\alpha !}(x'-y')^\alpha
\partial _{x'}^\alpha b(y',\xi ,\mu )v(y')\,dy'\d \xi \\
+\int e^{i(x'-y')\cdot \xi
'+ix_n\xi _n}\sum_{|\alpha |=M}\tfrac {M(x'-y')^\alpha}{\alpha
!}\int_0^1(1-h)^{M-1}
\partial _{x'}^\alpha b(y'+(x'-y')h,\xi ,\mu )v(y')\,dhdy'\d \xi \\
=\sum_{|\alpha |<M}\operatorname{OPK}(b_\alpha (y',\xi ,\mu
))v+\operatorname{OPK}(b_M(x',y',\xi ,\mu ))v;
\endmultline
$$
here we have replaced $(x'-y')^\alpha e^{i(x'-y')\cdot\xi '}$ by
$D_{\xi '}^\alpha e^{i(x'-y')\cdot\xi '}$ and integrated by parts
w.r.t.\ $\xi '$ (in the oscillatory integrals), and set
$$\aligned
b_\alpha (y',\xi ,\mu )&=\tfrac 1{\alpha !}\overline
D_{\xi '}^\alpha
\partial _{y'}^\alpha b(y',\xi ,\mu ),\\
b_M(x',y',\xi ,\mu )&=\sum_{|\alpha |=M}\tfrac M{\alpha !}\overline
D_{\xi '}^\alpha
\int_0^1(1-h)^{M-1}
\partial _{x'}^\alpha b(y'+(x'-y')h,\xi ,\mu )\,dh.
\endaligned\tag3.34$$

Now we compose with $\operatorname{OPT}(c_{jk}\bar{\hat\varphi
}_k)$ in front. For each term in the sum over $|\alpha |<M$, this
gives:$$
\gather
\operatorname{OPT}(c_{jk}\bar{\hat\varphi
}_k)\operatorname{OPK}(b_\alpha )=\operatorname{OP}'(s_\alpha
(x',y',\xi ',\mu )),\text{ where}\\
s_\alpha (x',y',\xi ',\mu )=\int c_{jk}\bar{\hat\varphi }_k\tfrac 1{\alpha
!}\overline
D_{\xi '}^\alpha
\partial _{y'}^\alpha h^+_{\xi _n}(q\hat\varphi _j)\,\d \xi _n
=\int c_{jk}\bar{\hat\varphi }_k\tfrac 1{\alpha !}\overline
D_{\xi '}^\alpha
\partial _{y'}^\alpha (q(y',\xi ,\mu )\hat\varphi _j)\,\d \xi _n.
\endgather
$$
In the last step it is used that $\overline D_{\xi '}^{\alpha
}\partial _{y'}^\alpha $ commutes with $h^+_{\xi _n}$ and that $\overline
D_{\xi '}^\alpha
\partial _{y'}^\alpha (q\hat\varphi _j)$ is a sum of terms with a
similar structure as $q\hat\varphi _j$, so that one can remove $h^+$
in the calculation of
the plus-integral, as in (B.5). More precisely,$$\gathered
\overline
D_{\xi '}^\alpha
\partial _{y'}^\alpha (q\hat\varphi _j)=\sum_{\gamma \le\alpha
}\tbinom\alpha \gamma \overline
D_{\xi '}^{\alpha-\gamma }
\partial _{y'}^\alpha q\overline D_{\xi '}^\gamma \hat\varphi _j,
\text{ where}\\
\overline
D_{\xi '}^{\alpha-\gamma }
\partial _{y'}^\alpha q=\sum_{m\le m'\le m+|\alpha -\gamma |}\frac
{r_{m'}(y',\xi )}{(p+\mu ^2)^{m'}},\quad
\overline D_{\xi '}^\gamma \hat\varphi _j=\sum_{|j'-j|\le |\gamma
|,j'\ge 0}r'_{j'}(\xi ')\hat\varphi _{j'};
\endgathered\tag3.35$$
here the $r_{m'}(y',\xi )$ are polynomials in $\xi $ of degree
$\le 2m'-J-2-|\alpha -\gamma |$
(vanishing if this number is $<0$) and
the $r'_{j'}(\xi ')$ are in $S^{-|\gamma |}$ (found from (B.2)). It
follows that$$
\overline
D_{\xi '}^\alpha
\partial _{y'}^\alpha (q(y',\xi ,\mu )\hat\varphi _j(\xi _n,\sigma ))=\sum
\frac {r_{\alpha
,m',j'}(y',\xi )}{(p+\mu ^2)^{m'}}\hat\varphi _{j'},\tag3.36
$$
where the sum is over $m\le m'\le m+|\alpha |$, $|j-j'|\le |\alpha |$,
$j'\ge 0$, and the $r_{\alpha ,m',j'}(y',\xi )$ are in
$S^{2m'-J-2-|\alpha |}(\Bbb R^{n-1},\Bbb R^n)$; they
are polynomials in $\xi
_n$ of degree $\le 2m'-J-2$, vanishing if this number is $<0$.
Each $r_{\alpha ,m',j'}(p+\mu ^2)^{-m'}$ is decomposed in
simple fractions with poles $\pm i\kappa ^\pm$ as in (A.7):
$$
\frac{r_{\alpha ,m',j'}(y',\xi )}{(p(y',\xi )+\mu
^2)^{m'}}= \sum_{1\le m''\le m'}\frac{r^+_{\alpha ,m',j',m''}(y',\xi ',\mu
)}{(\kappa ^+(y',\xi ',\mu )+ i\xi
_n)^{m''}}+\sum_{1\le m''\le m'}\frac{r^-_{\alpha ,m',j',m''}(y',\xi ',\mu
)}{(\kappa ^-(y',\xi ',\mu )- i\xi
_n)^{m''}};\tag3.37
$$
here the numerators are of order $m''-|\alpha |-J-2$, lying in
$S^{-|\alpha |,0,m''-J-2}$ (since the $\mu $-independent factors coming
from derivatives of $\hat\varphi _j$ have order down to $ -|\alpha |$),
and each $\lambda $-differentiation lowers the third upper index by $2$.

We can now apply Lemma B.4
to the compositions of $c_{jk}\bar{\hat \varphi }_k$ with (3.37)
(recalling from (B.1)  that the normalized Laguerre functions
have a factor $(2\sigma )^{\frac12}$). This gives $\psi $do symbols in
$(x',y')$-form on $\Bbb R^{n-1}$,
lying in $S^{\nu -|\alpha |,0,m''-J-2}\cdot
S^{1,0,-m''-1}\subset S^{\nu +1-|\alpha |,0,-J-3}$ if $j'\ne k$, and in
$S^{\nu -|\alpha |,0,-J-2}$ if $j'=k$.

Summing over the indices we find
that $s_\alpha (y',\xi ',\mu )\in
S^{\nu -|\alpha |,0,-J-2}$, and each $\lambda $-dif\-fe\-ren\-tia\-tion
lowers the
third upper index by $2$.

From the sequence of symbols $s_\alpha $, one constructs a symbol
$s\in S^{\nu ,0,-J-2}$ having the asymptotic expansion $s\sim
\sum_{\alpha \in\Bbb N^{n-1}}s_\alpha $ in the usual way.
There remains the question of how well this represents a symbol of the
operator given in (3.33). This will be dealt with by a consideration of
the
remainder after $M$ terms, using the formula in (3.34).

We shall use the strategy explained after (3.7).
In order to be able to distinguish between the various covariables,
we shall write here explicitly 
$\hat\varphi (\xi _n,[\xi '])$ and $\hat\varphi (\xi _n,[\eta '])$
instead of $\hat\varphi (\xi _n,\sigma )$.
Consider, for large $M$,
$$
\multline
\operatorname{OPT}(c_{jk}\bar{\hat\varphi
}_k)\operatorname{OPK}(b_M )v(x')\\
=\int_{\Bbb R^{4n-3}}
e^{i(x'-z')\cdot \xi '+i(z'-y')\cdot \eta '}c_{jk}(x',\xi
')\bar{\hat\varphi }_k(\xi
_n,[\xi ']) h^+_{\xi  _n}\bigl[\sum_{|\alpha |=M}\tfrac M{\alpha
!}\overline
D_{\eta  '}^\alpha\\
 \int_0^1(1-h)^{M-1}
\partial _{z'}^\alpha  q(y'+(z'-y')h,\eta ',\xi _n,\mu )
\hat\varphi _j(\xi
_n,[\eta '])\bigr]v(y')\,dh\d
\xi _ndy'\d\eta 'dz'\d\xi '\\
=\int_{\Bbb R^{4n-3}}
e^{i(x'-z')\cdot \xi '+i(z'-y')\cdot \eta '}c_{jk}\bar{\hat\varphi }_k(\xi
_n,[\xi '])
\sum_{|\alpha |=M}\tfrac M{\alpha !}\overline
D_{\eta  '}^\alpha
\int_0^1(1-h)^{M-1}\\
\times\partial _{z'}^\alpha  q(y'+(z'-y')h,\eta ',\xi _n,\mu )
\hat\varphi _j(\xi
_n,[\eta '])v(y')\,dh\d
\xi _ndy'\d\eta 'dz'\d\xi ';
\endmultline\tag3.38
$$
here we have again commuted $\overline D_{\eta  '}^{\alpha
}\partial _{z'}^\alpha $ and $h^+_{\xi _n}$, and removed the latter
because of the special structure as rational functions.
The kernel of this operator is
$$\aligned
&K_M(x',y ',\mu )=\int_{\Bbb R^{3n-2}}
e^{i(x'-z')\cdot\xi '+i(z'-y')\cdot \eta '}c_{jk}\bar{\hat\varphi }_k(\xi
_n,[\xi'])\sum_{|\alpha |=M}\tfrac M{\alpha !}\overline
D_{\eta  '}^\alpha\\
&\times \int_0^1(1-h)^{M-1}
\partial _{z'}^\alpha  q(y'+(z'-y')h,\eta ',\xi _n,\mu )
\hat\varphi _j(\xi
_n,[\eta '])\,dh\d\xi _nd\eta 'dz'\d\xi '.
\endaligned\tag3.39$$

Performing the differentiations $\overline D_{\eta  '}^{\alpha
}\partial _{z'}^\alpha $, one gets a sum of similar expressions, now
with order $M$ steps lower in the Poisson operator. Consider one term in
such
an expression, it is of the form (3.39) with $q$ replaced by a
fraction $r'(p+\mu ^2)^{-m'}$ as in (3.36), with $r'$ in
$S^{2m'-J-2-M}$:
$$\aligned
&\operatorname{OPT}(c_{jk}\bar{\hat\varphi}_k)
\operatorname{OPK}(b')v(x')
\\&\quad=\int_{\Bbb R^{4n-3}}
e^{i(x'-z')\cdot \xi '+i(z'-y')\cdot \eta '}c_{jk}\bar{\hat\varphi }_k(\xi
_n,[\xi '])\int_0^1(1-h)^{M-1}\\
&\quad\times
\frac{r'(y'+(z'-y')h,\eta ',\xi _n)}{(p(y'+(z'-y')h,\eta ',\xi _n)+\mu ^2)^{m'}}\hat\varphi _j(\xi
_n,[\eta '])v(y')\,dh\d
\xi _ndy'\d\eta 'dz'\d\xi '.
\endaligned\tag3.40
$$

Insert in (3.40) the expansion of the denominator from Lemma A.1:
$$
(p+\mu ^2)^{-m'}=\sum_{0\le i <L}c_{m',i }\mu ^{-2m'-2i }p^i
+\mu ^{-2m'-2L}p'_L.\tag 3.41
$$
The composed operators resulting from
the terms in the sum over $i $
are of the form $\mu ^{-2m'-2i }S_i $, where $S_i $ is a $\mu
$-independent
 $\psi $do on $\Bbb R^{n-1}$ of order
$\le\nu+2m'+2i -J-2-M$
(with $i <L$). These operators belong to the
calculus. Taking  $M$ large in comparison with $2L$ we can assure
that their kernels have any given finite degree of smoothness.

The contribution from $\mu ^{-2m'-2L}p'_L$ equals
$\mu ^{-2m'-2L} S'$, where $S'$ has a form as in (3.40) with $(p+\mu
^2)^{-m'}$ replaced by $p'_L$; it is a $\psi $do on $\Bbb R^{n-1}$
of order $\le\nu+2m'+2L -J-2-M$.
The estimates (A.2) imply that the entries satisfy
all symbol estimates {\it with uniform estimates in }$\mu \in\Gamma $
for $|\mu |\ge 1$.
Then $ S'$ maps $$
 S'\: H^{s}(\Bbb R^{n-1})\to H^{s-\nu-2m'-2L+J+2+M}(\Bbb R^{n-1})
, \text{ any }s\in\Bbb R,$$
with uniform estimates in $\mu \in\Gamma $ for $|\mu |\ge 1$, since the
norms
of such mappings are estimated in terms of the
symbol seminorms. Taking $M$ large in comparison with $2L$, we can
achieve
that the kernel of $S'$ is a smooth as we want, with uniform
bounds on the $(x',y')$-derivatives.

It follows that $$
S_{J,j,k,m}-\op'\bigl(s(x',y',\xi ',\mu
)\bigr)=S_{J,j,k,m}-\op'\bigl(\sum_{|\alpha
|<M}s_\alpha \bigr)-\op'\bigl(s-\sum_{|\alpha
|<M}s_\alpha \bigr)\tag3.42
$$
has a kernel expansion as in (3.7) (for $d=-J-2$), where $L$ and the
degree of smoothness can be as high as we want, since $M$ can be taken
arbitrarily large. 
Since differentiation in $\lambda $ lowers the $d$-index by two
steps, and $\partial _\lambda =-2\mu \partial _\mu $, the information
required for (3.8) is likewise found in this way.
This shows that $S_{J,j,k,m}$ is a $\psi $do with
symbol $\sim s$ in $S^{\nu ,0,-J-2}$. Its $\lambda
$-derivatives are treated in a similar way.

This shows the proposition in the case where the symbols are independent
of $x_n$. Now let $q_{-2-J}$ depend on $x_n$, and consider again
(3.32).
By a Taylor expansion in $x_n$, $q_{-2-J}$ is of the form$$
\multline
q_{-2-J}(x,\xi ,\mu )=\sum_{0\le r<M}\tfrac 1{r!}x_n^r\partial
_{x_n}^rq_{-2-J}(x',0,\xi ,\mu )\\
+\tfrac1{(M-1)!}{x_n^M}\int_0^1(1-h)^{M-1}
\partial _{x_n}^{M}q_{-2-J}(x',hx_n,\xi ,\mu )\,dh.
\endmultline\tag 3.43$$
For the terms in the sum over $r$,
we observe that the factor $x_n^r$
goes together with the trace operator to the left
$\operatorname{OPT}(c_{jk}(x',\xi ')\bar{\hat\varphi }_k(\xi
_n,\sigma ))$ and replaces its symbol by $c_{jk} D_{\xi
_n}^r\bar{\hat\varphi
}_k(\xi _n,\sigma )$ [G2, (2.4.14)]. This lowers the degree --- but not
the $d$-index --- by $r$ steps; more precisely it gives $c_{jk}[\xi
']^{-r}$ times a linear combination of adjacent Laguerre functions,
cf\. (B.2). The resulting composition is of the type we
have dealt with above, giving a symbol in $S^{\nu -r,0,-2-J}$.

The last term in (3.43) is a linear combination of terms of the form
$$
x_n^M\int_0^1(1-h)^{M-1}\partial
_{x_n}^{M}\bigl(\frac{r(x',hx_n,\xi )}{(p(x',hx_n,\xi )+\mu
^2)^m}\bigr)\,dh,
$$
with $r$ polynomial in $\xi $ of order $2m-2-J\ge 0$;
consider one such term. Insert the expansion (3.41) of
$(p+\mu ^2)^{-m}$.
The composed operators resulting from
the terms in the sum over $i $
are of the form $\mu ^{-2m-2i }S_i
$
with $S_i $ of order $\le \nu +2m+2i -2-J-M$
(with $i <L$). The operators belong to the
calculus, with more smoothness of the kernels, the larger $M$ is
taken in comparison with $2L$.
Finally, the contribution from $p'_L$ is of the form $\mu
^{-2m-2L}\widetilde S$, with
$$
\gathered
\widetilde S=\operatorname{OPT}(c_{jk}D_{\xi _n}^M
\bar{\hat\varphi}_k(\xi _n,\sigma))
\operatorname{OP}( \tilde q_{L}(x,\xi ,\mu))_+
\operatorname{OPK}(\hat\varphi _j(\xi_n,\sigma )),\\
\tilde q_{L}(x,\xi ,\mu )=
\tfrac1{(M-1)!}\int_0^1(1-h)^{M-1}\partial
_{x_n}^{M}\bigl({r(x',hx_n,\xi )}p'_{L}(x',hx_n,\xi ,\mu
)\bigr)\,dh.\endgathered\tag 3.44
$$
The estimates (A.2) imply that for all indices:$$
|D_{x}^\beta D_\xi ^\alpha \tilde q_{L}(x,\xi ,\mu )|\le
C'_{\alpha ,\beta }\ang\xi ^{2m-2-J+2L-|\alpha |},\tag 3.45
$$
with $C'_{\alpha ,\beta }$ independent of $\mu \in\Gamma $ for $|\mu |\ge
1$.
Then the resulting operator $\widetilde S$ maps $H^{s}(\Bbb R^{n-1})$
to $ H^{s-\nu-2m+2+J-2L+M+1}(\Bbb R^{n-1})$
for any $s\in\Bbb R$,
with uniform estimates in $\mu \in\Gamma $ for $|\mu |\ge 1$. Taking
$L$ large,
and $M$ large
in comparison with $2L$, we can
obtain that the resulting kernel is a smooth as we want, with
coefficient $\mu ^{-2m-2L}$ of as low order as we want.

As in the preceding considerations, this allows us to conclude that
the whole composed operator belongs to the calculus, with symbol
properties as
asserted.
\qed
\enddemo

\example {Remark 3.9} In the above proof, one can moreover keep track
of the size of symbol
seminorms in their dependence on $j$ and $k$. It follows from Lemma B.4
and the results on $\xi '$-derivatives and $\xi _n$-derivatives we
used, that for each $J$, the symbol seminorms resulting from calculations
with Laguerre funtions are only polynomially increasing in $j$ and
$k$. Thus if the symbol seminorms of the $c_{jk}$ are rapidly decreasing
in $j$ and $k$, as they are when the $c_{jk}$ come from an s.g.o.\ as
in Lemma B.6, summations in $j$ and $k$ of operators $S_{J,j,k}$ will
converge in the relevant symbol seminorms.
\endexample

Note that for $J\ge 1$, the $(N-1)$'st derivatives of the symbols of the
operators $S_{J,j,k}$ lie in $S^{\nu ,0,-2N-1}\subset S^{\nu
-2N-1,0,0}\cap S^{\nu ,-2N-1,0}$ (cf\. (3.3));
hence in their diagonal kernel expansions as in
(3.18), the sum over $k$ starts with $k\ge 1$, so that they
contribute only locally to the coefficient of $\mu ^{-2N}$. Thus only
the terms with $J=0$,
$$
\aligned
\tfrac{\partial _\lambda ^{N-1}}{(N-1)!}S_{0,j,k}&=
\tfrac{\partial _\lambda ^{N-1}}{(N-1)!}C_{jk}\Phi
^*_k\op((p_{1,2}+\mu ^2)^{-1})_+\Phi _j\\
&=C_{jk}\Phi
^*_k\op((p_{1,2}+\mu ^2)^{-N})_+\Phi _j
\endaligned \tag3.46
$$
can contribute nonlocally to the coefficient of $\mu ^{-2N}$.
We can use the analysis in the proof of Proposition 3.8 to
distinguish the contributing part still further:

\proclaim{Proposition 3.10} We have for each $N\ge 1$ that$$
C_{jk}\Phi
^*_k\op((p_{1,2}+\mu ^2)^{-N})_+\Phi _j
=\op'\big(\delta _{jk}c_{jk}(x',\xi ')\alpha ^{(N)}(y',\xi ',\mu
)\big)+\widetilde S_{0,j,k,N},\tag3.47
$$
where $\alpha ^{(N)}$ is as defined in {\rm (A.12)} and $\widetilde
S_{0,j,k,N}$ has symbol in $S^{\nu +1,0,-2N-1}$.
\endproclaim

\demo{Proof} First let $q$ be independent of $x_n$, and consider how
$q_{-2}=(p+\mu ^2)^{-1}$ (with $p=p_{1,2}$) enters in the proof of
Proposition 3.8, for $N=1$. Whenever we apply a differentiation
$\partial _{x'}^\alpha $ with $|\alpha |>0$ to $q$, it produces a sum
of rational functions
 with denominators $(p+\mu ^2)^{-m}$ with $m\ge 2$ and $\mu
$-independent numerators, so only the undifferentiated term $q_{-2}$
retains the power $-1$. Therefore the resulting terms in (3.35) with
$|\alpha |>0$ have denominators $(p+\mu ^2)^{-m'}$ with $m'\ge 2$.
Then when we decompose in simple fractions and apply Lemma B.4, we
get symbols lying in $S^{\nu +1,0,-3}$. 
There remains the term$$
 c_{jk}(x',\xi ')\int\bar{\hat\varphi }_k(\xi _n,\sigma )(p(y',\xi
',\xi _n)+\mu ^2)^{-1}\hat\varphi _j(\xi _n,\sigma )\,\d\xi _n.
$$
Here Lemma B.4 shows that it gives a symbol in $S^{\nu +1,0,-3}$ when
$j\ne k$, whereas it gives $c_{jk}(x',\xi ')\alpha ^{(1)}(y',\xi
',\mu )\in S^{\nu ,0,-2}$ when $j= k$ (note that the norming factors
$(2\sigma )^\frac12$ in the Laguerre functions eliminate the division
by $2\sigma $ in (B.8)).

This shows the claim for $N=1$, and in view of the information on
$\lambda $-derivatives, it follows for general $N$.

When $q$ depends on $x_n$, we consider the Taylor expansion (3.43).
Again, any derivative $\partial _{x_n}^r$ with $r>0$ produces
rational functions with $p+\mu ^2$ in powers $\ge 2$ in the
denominators. So we find that the Taylor terms, except for the first
one, lead to symbols with lower third index. The first term is dealt
with above.
\qed
\enddemo

\proclaim{Theorem 3.11} Consider a localized situation, where $G$ is given
in the form  {\rm (3.30)} and the symbols $c_{jk}(x',\xi ')$ have compact
$x'$-support. Define the diagonal sum$$
C=\sum_{j\in\Bbb N}C_{jj}, \text{ with symbol }c(x',\xi
')=\sum_{j\in\Bbb N}c_{jj}(x',\xi ').\tag 3.48
$$
Then
$$\Tr _{\Bbb R^n_+}(GQ^N_{\lambda ,+})=\Tr_{\Bbb
R^{n-1}}\op'(c(x',\xi ')\alpha ^{(N)}(x',\xi ',\mu ))+\Tr_{\Bbb
R^{n-1}}\widetilde S',\tag 3.49$$
where $\widetilde S'$
has symbol in $S^{\nu +1,0,-2N-1}$. It follows that for $N>(\nu
+n-1)/2$, the trace has an expansion (with all $c'_k=0$ if $\nu
\notin\Bbb Z$): 
$$
\Tr_{\rnp}(GQ^N_{\lambda ,+})\sim \sum_{j\ge 1}\! c_{j}(-\lambda)
^{\frac{n+\nu - j}{2}-N}
+\sum_{k\ge
0}( c'_k\log (-\lambda )+ c''_k)(-\lambda )
^{-\frac k2-N},
\tag 3.50$$
 where  $$
c''_0=\int \slint \tr c(x',\xi ')\,\d\xi 'dx'+\text{
local terms}\tag3.51
$$
(with vanishing local terms if $\nu \in\Bbb
R\setminus \Bbb Z $ or $\nu <1-n$).
\endproclaim

\demo{Proof}
For each $j,k$, we conclude from
Proposition 3.10 that
$$
\aligned
\Tr _{\Bbb R^n_+}(\Phi _jC_{jk}\Phi ^*_kQ^N_{\lambda ,+})&=
\Tr _{\Bbb R^{n-1}}(C_{jk}\Phi ^*_kQ^N_{\lambda ,+}\Phi _{j})
\\&=\Tr_{\Bbb R^{n-1}}\op'(\delta _{jk}c_{jk}(x',\xi ')\alpha
^{(N)}(y',\xi ',\mu
))+\Tr_{\Bbb R^{n-1}}\widetilde S'_{0,j,k,N}
\\ &=\Tr_{\Bbb
R^{n-1}}\op'(\delta _{jk}c_{jk}(x',\xi ')\alpha ^{(N)}(x',\xi ',\mu
))+\Tr_{\Bbb R^{n-1}}\widetilde S'_{0,j,k,N}.
\endaligned
$$
with pseudodifferential operators $\tilde S'_{0,j,k,N}$ having symbols
in $S^{\nu+1,0,-2N-1}$.
For the last equality
it is used that the trace is calulated from
the diagonal values of the kernel, where $y'$ is taken equal to $x'$.
The formula (3.49) follows by summation over
$j$ and $k\in\Bbb N$, using that $(C_{jk})_{j,k\in\Bbb N}$ is rapidly
decreasing in
$j$ and $k$ and that the symbol seminorms of the composed operators
are polynomially controlled in $j$ and $k$, as indicated in Remark 3.9.

Since $c_{jj}\alpha^{(N)}\in S^{\nu-2N,0,0}\cap S^{\nu,-2N,0}$,
the trace expansion follows from (3.5) with the statement (3.51)
as in the proof of Theorem 3.3.
(We have replaced $j$ by $j+1$ in
(3.18) and inserted $\mu =(-\lambda )^{\frac12}$, to facilitate
comparison with (1.10).)
\qed\enddemo

We have as an immediate corollary (using Proposition B.3 for
$GG^{(N)}_\lambda $):

\proclaim{Theorem 3.12} Let $G$ be a singular Green operator on $X$
of order $\nu \in\Bbb R$ and class $0$.
Consider a localization $\underline G$ to $\crnp$ as described after
Lemma {\rm 3.1}, write$$
\underline G=\sum_{j,k\in\Bbb N}\Phi _j\underline C_{jk}\Phi ^*
_k,\quad \underline C= \sum_{j\in\Bbb N}\underline C_{jj},\tag3.52
$$
with $\underline C=\op'(c(x',\xi '))$, and carry $\underline C$ back
to an operator $C$ on $X'$ by use of the considered
coordinate mappings.
Then
$$
C_0(G,P_{1,\operatorname{D}})=\int\slint
\tr c(x',\xi ')\,\d\xi ' dx'
+\text{ local terms,}\tag 3.53
$$
with
vanishing local terms if $\nu <-n+1$ or $\nu \in\Bbb R\setminus \Bbb
Z$. In particular, we have for the localized operators as well as for
the operators on $X,X'$:
$$
\aligned
C_0(\underline G,\underline P_{1,\operatorname{D}})&=C_0(\underline
C,\underline S)+\text{ local terms,}\\
C_0(G,P_{1,\operatorname{D}})&=C_0(C,S)+\text{ local terms,}
\endaligned\tag 3.54
$$
for any auxiliary elliptic second-order operator $S$ on $X'$ with no
principal symbol eigenvalues on $\Bbb R_-$.

\endproclaim

This can be used to show the commutation property for
$C_0(G,P_{1,\operatorname{D}})$:

\proclaim{Theorem 3.13} When $G$ and $G'$ are
singular Green operators of orders $\nu$ resp\. $\nu'$ and class 0,
then $C_0([G,G'],P_{1,\operatorname{D}})$ is locally determined; it
vanishes if $\nu +\nu '\notin\Bbb Z$ or $\nu +\nu '<1-n$.
\endproclaim

 \demo{Proof}
Expand the localized versions in series with Laguerre operators,
$$\underline G=\sum_{j,k\in\Bbb N}\Phi _j\underline C_{jk}\Phi ^*
_k,\quad
\underline G'=\sum_{j,k\in\Bbb N}\Phi _j\underline C'_{jk}\Phi ^*_k.
$$
In view of (B.13),$$\aligned
[\underline G,\underline G']&=\sum_{j,k,l,m}\Phi _j \underline C_{jk} \Phi
_k^*\Phi _l
\underline C'_{lm} \Phi _m^* - \sum_{j,k,l,m}\Phi _j \underline C'_{jk}
\Phi
_k^*\Phi _l \underline C_{lm} \Phi _m^*\\
&=\sum_{j,k,m}\Phi _j (\underline C_{jk} \underline
C'_{km}-\underline C'_{jk} \underline C_{km})
 \Phi _m ^* ,
\endaligned$$
and the  diagonal sum associated with this operator as in (3.52) is:
$$\sum_{j,k} (\underline C_{jk} \underline C'_{kj}-\underline C'_{jk}
\underline C_{kj})=\sum_{j,k} \underline C_{jk} \underline
C'_{kj}-\sum_{j,k}\underline C'_{kj}
\underline C_{jk}=\sum_{j,k} [\underline C_{jk}, \underline C'_{kj}].$$
 This series converges in the relevant symbol seminorms, in view of the
rapid decrease of the $\underline C_{jk}$ and $\underline C'_{jk}$ for
$j,k\to\infty $. Then
it follows from Theorem 3.12 that
$$
C_0([\underline G,\underline G'],\underline
P_{1,\operatorname{D}})=C_0(\sum_{j,k}[\underline C_{jk},
\underline C'_{kj}],\underline S)+ \text{ local terms }.
$$
Now we use the fact known for closed manifolds (cf\. e.g. \cite{G4})
that each $C_0([\underline C_{jk},\underline C'_{kj}],\underline S)$
is local, defined from specific homogeneous terms in the symbols;
then so is the sum. In particular, they vanish in the cases  $\nu
+\nu '\notin\Bbb Z$ or $\nu +\nu '<1-n$; then so does the sum. The
resulting statements carry back to the manifold situation.
Hence $C_0([G,G'], P_{1,\operatorname{D}})$ is
local.
\qed\enddemo

It may be remarked that when $G=\sum _{j,k\in\Bbb N}\Phi
_jC_{jk}\Phi ^*_k$ on $\crnp$ with $C_{jk}=\op'(c_{jk}(x',\xi '))$, then
the symbol of $G$ has the same principal
part as $\sum_{j,k\in\Bbb N}\hat\varphi _j(\xi _n,\sigma
)c_{jk}(x',\xi ')\bar{\hat\varphi }_k(\eta _n,\sigma )$, but need not
equal the full symbol. Likewise, the normal trace $\tr_n G$ has the
same principal symbol as the diagonal sum
$\sum_{j\in\Bbb N}C_{jj}$, but not always the same full symbol.
In any case,
the above considerations show that they give the same
contribution to $C_0(G,P_{1,\operatorname{D}})$ (modulo local terms,
if $\nu \in\Bbb Z$, $\nu \ge -n$).

This ends the proof that $C_0(G,P_{1,\operatorname{D}})$ is a
quasi-trace on s.g.o.s. We can also describe cases where it is a
trace; more about that below.

The local terms, we have talked about so far, are defined from the
symbols on $X$. But in fact, it is only the behavior in an
arbitrarily small neighborhood of $X'$ that enters:

Assume that a normal coordinate $x_n$ has been chosen such that a suitable
neighborhood $X_1$ of $X'$ in $X$ is represented by the product
$X_1=X'\times [0,1]$, with $x\in X_1$ written as $x=(x',x_n)$,
$x'\in X'$ and $x_n\in [0,1]$. We can also assume that $E|_{X_1}$ is
the lifting of $E|_{X'}$. Let $\chi (t)$ be a $C^\infty $
function on $\crp$ that is 1 for $t\le 1$ and $0$ for $t\ge 2$, and,
for $\varepsilon <\frac12$,
let $$
\chi _\varepsilon (x)=\cases \chi (x_n/\varepsilon )\text{ for }x\in
X_1,\\
0\text{ for }x\in X\setminus X_1.\endcases\tag3.55
$$
Then we can write $$
G=G_{\operatorname{b}}+G_{\operatorname{i}}; \quad
G_{\operatorname{b}}=\chi _{\varepsilon /4}G\chi _{\varepsilon
/4},\tag3.56
$$
here $G_{\operatorname{b}}$ is supported in $X_{\varepsilon
/2}=X'\times [0,\varepsilon /2]$,
and $G_{\operatorname{i}}$ is of order $-\infty $ (in
particular it is trace-class).
The auxiliary local coordinates can be chosen such that the variable $x_n$
is
preserved for the patches intersecting with $X_1$; then $\tr_n
G_{\operatorname{b}}$ has a meaning as a $\psi $do on $X'$.

\proclaim{Corollary 3.14}

{\rm (i)} With the decomposition defined above,
$$
C_0(G,P_{1,\operatorname{D}})=C_0(G_{\operatorname{b}},P_{1,D})+\Tr_X
G_{\operatorname{i}}+ \text{ local terms},\tag 3.57
$$
where the local terms depend on the symbol of $P_1$ and the first
$\nu +n$ symbols of $G$ on $X_\varepsilon $ only.

{\rm (ii)} If $\nu \in \Bbb R\setminus \Bbb Z$ or $\nu <-n+1$, then $$
C_0(G,P_{1,\operatorname{D}})=\operatorname{TR}_{X'}(\tr_n
G_{\operatorname{b}})+\Tr_X
G_{\operatorname{i}}.\tag 3.58
$$
For such values of $\nu $,
$C_0(G,P_{1,\operatorname{D}})$ is a canonical trace on the
operators $G$:
$$
C_0(G,P_{1,\operatorname{D}})=\TR G,\tag3.59
$$ in the sense that it satisfies:
\roster
\item"(1)" $C_0(G, P_{1,\operatorname{D}})$ is independent of $P_1$
if $\nu \notin\Bbb Z$ or $\nu <1-n$;
\item"(2)" $C_0([G,G'],P_{1,\operatorname{D}})=0$ if $\nu +\nu
'\notin\Bbb Z$ or $\nu +\nu '<1-n$.
\endroster

{\rm (iii)} If $\nu < 1-n$, $$
\TR G=\Tr_X G.\tag 3.60
$$
\endproclaim

\demo{Proof}  For (i), note that since
$\chi _{\varepsilon /4}=\chi _{\varepsilon /2}\chi _{\varepsilon /4}
$
$$
\Tr(G_{\operatorname{b}}R_\lambda ^N)=\Tr(\chi _{\varepsilon /2}\chi
_{\varepsilon /4}G\chi _{\varepsilon /2}R_\lambda ^N)=
\Tr(\chi _{\varepsilon /4}G\chi _{\varepsilon
/2}R_\lambda ^N\chi _{\varepsilon /2}),
$$
which shows that only the behavior of $ R_\lambda $ on
$X_\varepsilon $ and the behavior of $G$ on $X_{\varepsilon /2}$
enter in the calculation of $C_0(G_{\operatorname{b}},
P_{1,\operatorname{D}})$.

As noted in the preceding theorems, the formulas are exact without
unspecified local terms, when $\nu $ is as in (ii). Moreover,
$C_0(G_{\operatorname{b}},
P_{1,\operatorname{D}})=C_0(\tr_nG_{\operatorname{b}},S)$
then by (3.29), for the allowed auxiliary elliptic operators $S$.
For the manifold $X'$ without boundary we have
from \cite{KV}, \cite{L}, \cite{G4} that $C_0(\tr_n
G_{\operatorname{b}},S)=\TR
(\tr_n G_{\operatorname{b}})$ for such $\nu $. This shows (3.58).
The statements in (1) and (2) have been shown further above. 

For (iii), we note that when $\nu <1-n$, $G_{\operatorname{b}}$
is trace-class with continuous kernel, and $$
\Tr G_{\operatorname{b}}=\Tr _{X_\varepsilon
}G_{\operatorname{b}}=\int_{X'}\int_0^\varepsilon
K(G_{\operatorname{b}}, x,x)\,dx_ndx'=\Tr_{X'}(\tr_n
G_{\operatorname{b}}).\tag3.61$$
Then $C_0( G,P_{1,\operatorname{D}})=\Tr G_{\operatorname{b}}+\Tr
G_{\operatorname{i}}=\Tr G$.
\qed
\enddemo

If $G_{\operatorname{b}}$ is written in the local coordinates in the
composition form (3.52), one can show similar formulas with $C$
instead of $\tr_nG_{\operatorname{b}}$.

When $\nu $ is integer, it may be so that $\tr_n
G_{\operatorname{b}}$ or $C$ has a parity that fits with the
dimension $n-1$, cf\. (1.4)--(1.5)ff. One could then try to show a
canonical trace property of $C_0(G,P_{1,\operatorname{D}})$ as in
the corresponding situation for operators given directly on the
closed manifold $X'$ as
in \cite{G4}. This is further supported by the fact that
since $\kappa ^+(x',-\xi ',\mu )=\kappa ^{-}(x',\xi ',\mu
)$,  $\alpha ^{(1)}$ is even-even, cf\. (A.12). However, there are a
lot of terms that give local contributions; not only lower order
terms in compositions but notably all the  $c_{lm}(x',\xi ')$ with
$l\ne m$ in the symbol $\sum_{l,m\in\Bbb N}c_{lm}\hat\varphi
_l\bar{\hat\varphi }_m$ of $G$; cf\. Lemma B.4: For $l>m$, the
resulting symbols involve $\kappa ^+$, for $l<m$ $\kappa ^-$, on the
principal level.
A high degree of symmetry seems required, or we can restrict the
auxiliary operators to a nicer type:

Consider $G$
as defined on the cylinder $X'\times \crp$ (after being cut down to a
neighborhood $X_\varepsilon $ of $X'$). Referring to this fixed
choice of normal coordinate,
we assume that the symbol of $P_1$ is not only even-even with respect to
$\xi $ (as any differential operator is), but it is so with respect
to $\xi '$. This holds when there are no terms with $D_{x_n}$ times a
first-order operator in $x'$, for example when$$
P_1=D_{x_n}^2+P'_1,\tag3.62
$$
where $P'_1$ is a second-order elliptic differential operator on $X'$.
We can take $P'_1$ with positive principal symbol, to assure strong
ellipticity of $P_1$ (a little more generality could be allowed).

Assume that $G=\opg(\sum_{j,k\in\Bbb N}c_{jk}\hat\varphi
_j\bar{\hat\varphi }_k)$ (expressed in local coordinates on $X'$,
using $x_n$ as the normal coordinate), with the $c_{jk}(x',\xi ')$
all being even-even or all being even-odd.  Then when
$GQ^N_{\lambda ,+}$ and $GG^{(N)}_\lambda $ are calculated, we find
that $$ \tr_n GQ^N_{\lambda ,+}\text{ and }\tr_n GG^{(N)}_\lambda $$
 are even-even of order $\nu -2N$, resp\. even-odd of order $\nu
-2N$. An application of the proof of \cite{GS1, Th\. 2.1} to these
two $\psi $do's on $X'$ gives that in the even-even case, the local
contributions to the coefficient of $\mu ^0$ and to $\log \mu $
vanish if $n-1$ is odd, and in the even-odd case, the local
contributions to the coefficient of $\mu ^0$ and to $\log \mu $
vanish if $n-1$ is even. So indeed, restricting the auxiliary
operators to those of the form (3.62), we have a canonical trace in
suitable parity cases. We have shown:

\proclaim{Theorem 3.15} Let $G$ be given on the cylinder $X'\times \crp$
with points $(x',x_n)$
and let $P_1$, in addition to being strongly elliptic with scalar
principal symbol, be of the form {\rm (3.62)}.
Let  $G$ be of order $\nu \in\Bbb Z$, with symbol  $\sum_{j,k\in\Bbb
N}c_{jk}(x',\xi ')\hat\varphi
_j(\xi _n,\sigma )\bar{\hat\varphi }_k(\eta _n,\sigma )$ (for a choice
of local
coordinates on $X'$, with $x_n$ as normal coordinate). Then $$
C_0(G,P_{1,\operatorname{D}})=\TR (\tr_n G)\tag3.63
$$
holds in the following cases:

{\rm (i)} All the $c_{jk}$ are even-even and $n$ is even.

{\rm (ii)} All the $c_{jk}$ are even-odd and $n$ is odd.
\endproclaim

In this sense, we may say that s.g.o.s have the canonical trace
$$
\TR G=\TR_{X'}(\tr_n G),\tag3.64
$$in
the above cases (i) and (ii).
One can also show that $C_0([G,G'], P_{1,\operatorname{D}})=0$ when,
in the case $n$ even, the $c_{jk}$ and $c'_{jk}$ are either all
even-even or all even-odd; in the case $n$ odd it holds when all the
$c_{jk}$ are even-even and all the $c'_{jk}$ are even-odd.

If $G$ is given in the form (3.52), the
statements hold with
 $\tr_n G$
replaced by $C=\sum_{j\in\Bbb N}C_{jj}$.

One could also introduce a parity concept
directly for s.g.o.\ symbols $g$, but it involves some unpractical
shifts: Recall
that in the polyhomogeneous expansion $g\sim\sum_{l\ge 0}g_{\nu -1-l}$,
the
homogeneity  degree of the $l$'th term is $\nu -1-l$, one step lower
than the order. Since the $\varphi _j$
depend on $\xi '$
through $\sigma (\xi ')=[\xi ']$, the parity statements on the $c_{jk}$
correspond to the opposite parity in $\xi '$ for $g$, relative to the
degrees. --- Observe moreover that the parity property depends on the
choice
of normal coordinate (cf\. the general transformation rule \cite{G2,
(2.4.62)ff.}).

For an example where the result may be of interest, let us mention that
the
singular Green part of the solution operator for the Dirichlet or
Neumann problem for a strongly elliptic differential operator of the form
(3.62) has even-even parity for
the coefficients $c_{jk}$ in the Laguerre expansion of the symbol.
(This is seen by calculations related to those in Remark 4.2 below, for
fixed
$\lambda $.)
Then it has a canonical trace if $n$ is even.

If we merely have that $\tr_n g$ has a parity that fits the dimension
$n-1$ (instead of all $c_{jk}$ having it), then there holds at least
that $\TR (\tr_n G)=C_0(\tr_n G,S)$
is well-defined independently of $S$ when $n$ is
even, resp\. odd, by the result for closed manifolds (\cite{KV},
\cite{G4}) applied to $X'$. This number can then be viewed as the
nonlocal part of $C_0(G,P_{1,\operatorname{D}})$, in an explicit
way.

\example{Remark 3.16} The application of \cite{GS1, pf. of Th. 2.1}
to $\tr_n GR^N_\lambda $ shows even more, namely that in the parity cases (i)
and (ii) in Theorem 3.15,  all the coefficients $\tilde c_j$ 
with $\nu +n-1-j$ even and all the log-coefficients $\tilde c'_k$ with
$k$ even
vanish in (1.10) (cf\. also \cite{G4, pf.\ of Th\. 1.3}). This holds also for
the pointwise kernel 
expansions as in (3.18).
\endexample

Finally, let us mention another useful application of the method
leading to Theorem 3.13:

\proclaim{Corollary 3.17} When $K$ is a Poisson operator of order $\nu
$ and $T$ is a trace operator of class $0$ and order $\nu '$, then,
for auxiliary operators $P_{1,\operatorname{D}}$ and $S$ as above,
$$
C_0(KT,P_{1,\operatorname{D}})=C_0(TK,S) +\text{ local terms};\tag 3.65
$$
the local terms vanishing if $\nu +\nu '<1-n$ or $\notin\Bbb Z$.
\endproclaim

\demo{Proof} Applying a partition of unity,
we can reduce to the case where the operators are given on
$\rnp$ (with compact $x'$-support).
Here we write (cf\. Lemma B.7)
$$
K=\sum_{j\in\Bbb N}\Phi _jC_j,\quad T=\sum_{k\in\Bbb N}C'_k\Phi _k^*,
$$
with rapidly decreasing sequences $(C_j)_{j\in\Bbb N}$ and
$(C'_k)_{k\in\Bbb N}$ of $\psi $do's on $\Bbb R^{n-1}$ of order $\nu
-\frac12$ resp\. $\nu '+\frac12$ (and with compact $x'$-support). By
(3.54) and (B.13),$$\aligned
C_0(KT, P_{1,\operatorname{D}})&=C_0(\sum_{j,k\in \Bbb N}\Phi
_jC_jC'_k\Phi  _k^*, P_{1,\operatorname{D}})\\
&=C_0(\sum_{j}C_jC'_j,S)+\text{ local terms};\\
C_0(TK,S)&=C_0(\sum_{j,k\in\Bbb N}C'_k\Phi _k^*\Phi _jC_j,S)
=C_0(\sum_{j}C'_jC_j,S).
\endaligned$$
The latter equals $C_0(\sum_{j}C_jC'_j,S)$ modulo local terms, by the
commutativity property for $\psi $do's on $\Bbb R^{n-1}$ and the
rapid decrease of the $C_j$ and $C'_j$. This shows
the identity (3.65), and the vanishing statement follows in the usual
way.
\qed
\enddemo

\head 4. The integer order case with nonvanishing $\psi $do part
\endhead

Consider now the case where $A$ equals $P_++G$, of order $\nu \in\Bbb
Z$. Here
$C_0(P_++G,P_{1,\operatorname{D}})=C_0(P_+,P_{1,\operatorname{D}})
+C_0(G,P_{1,\operatorname{D}})$,
where $C_0(G,P_{1,\operatorname{D}})$ has already been analyzed. For
the determination of $C_0(P_+,P_{1,\operatorname{D}})$, we rewrite as
follows:
$$\aligned
\Tr(P_+R_\lambda ^N)&=\Tr(P_+Q_{\lambda ,+} ^N)+\Tr(P_+G_\lambda
^{(N)})\\
&=\Tr((PQ^N_{\lambda })_+)-\Tr(G^+(P)G^-(Q^N_{\lambda }))+\Tr(P_+G_\lambda
^{(N)}).
\endaligned\tag4.1$$
It follows from the results of \cite{GSc} recalled in Appendix B
that the last two terms have expansions as in Proposition B.3,
contributing only locally to $C_0(P_+,P_{1,\operatorname{D}})$. For
the first term, we find the desired information by considering
$PQ^N_\lambda $ on $\widetilde X$.
As in Section 3 (see the details after Lemma 3.1), we can use local
coordinates and a subordinate
parition of unity, now covering all of $\widetilde X$ with open sets
$U_i$, such that subsets of
$X$ resp\. $\widetilde X\setminus X$ are mapped to subsets
of $\crnp$ resp\.
$\rnm$, the intersections of the $U_i$ with $X'$ being mapped into
$\partial \crnp$. Again we can replace the images $V_{j(i_1,i_2)}$ by
sets $V'_{i_1,i_2}$ with positive distance from one another (in
the $x'$-direction), so
that we can refer to one localized operator.
By \cite{GS1, Th\. 2.1 and 2.7}
and the more detailed information in \cite{G4, Th\. 1.3 ff.}, the
kernel of $PQ^N_\lambda $ in the localized
situation
has a diagonal expansion:$$
K( PQ^N_\lambda ,x,x)
\sim
\sum_{ j\in \Bbb N  } \tilde c_{ j}(x)(-\lambda )  ^{\frac{\nu +n -j}2-N}+
\sum_{k\in \Bbb N}\bigl( \tilde c'_{ k}(x)\log (-\lambda )  +\tilde c''_{
k}(x)\bigr)(-\lambda ) ^{-k-N},
\tag4.2
$$
where $$
c_{\nu +n}(x)+c''_0(x)=\slint p(x,\xi )\,\d\xi +\text{ local
terms}.\tag4.3
$$
From the sets mapped into $\crnp$, we find
$\Tr_{X}((PQ^N_\lambda )_+)$ by integrating over ${\Bbb R^n_+}$. Then
(4.2) implies
$$
\Tr_X((PQ^N_\lambda )_+)\sim \sum_{ j\in \Bbb N  } \tilde c_{ j,+}(-\lambda
)  ^{\frac{\nu +n -j}2-N}+
\sum_{k\in \Bbb N}\bigl( \tilde c'_{ k,+}\log (-\lambda )  +\tilde c''_{
k,+}\bigr)(-\lambda ) ^{-k-N},
\tag4.4
$$
where the coefficients are obtained from the coefficients in (4.2) by
integration in $x$. In particular, in view of (4.3),$$
c_{\nu +n,+}+c''_{0,+}=\int_{\Bbb R^n_+} \slint  \tr p(x,\xi )\,\d\xi dx
+\text{ local terms.}\tag4.5
$$
Here $c_{\nu +n,+}$ is also locally determined, so we can conclude:

\proclaim{Theorem 4.1} Let $P$ be of order $\nu \in\Bbb Z$, having
the transmission property at $X'$. In the localized situation,
$$
C_0(P_+,P_{1,\operatorname{D}})=\int_{\Bbb R^n_+} \slint \tr p(x,\xi
)\,\d\xi
dx+\text{ local terms.}\tag4.6
$$
When
$P_1$ and $P_2$ are two choices of auxiliary operator, then
$C_0(P_+,P_{1,\operatorname{D}})-C_0(P_+,P_{2,\operatorname{D}})$ is
locally determined.
\endproclaim

\demo{Proof} It only remains to establish the last statement: It
holds, since the global term
$\int_{\Bbb R^n_+} \tslint \tr p(x,\xi )\,\d\xi dx$
cancels out in the calculation of the
difference in local coordinates, leaving only locally determined
terms.\qed
\enddemo

\example{Remark 4.2}
As usual, we can ask for cases where the local terms vanish. They do
so for $\nu <-n$, where $$
C_0(P_+,P_{1,\operatorname{D}})=\Tr P_+
$$
is easily seen. For $\nu \ge -n$, one can analyze cases with parity
properties, but they will in general not have vanishing local
contributions, as the following considerations show:
Let $P_+=I$, the simplest possible choice, and let
$P_1$ be of the simple form (3.62) with $P'_1$ selfadjoint positive. We
have
to expand$$
\Tr(IR^N_\lambda )=\Tr (Q^N_{\lambda ,+})+ \Tr (G^{(N)}_\lambda ).
$$
Since $Q_\lambda $ is even-even,$$
\Tr (Q^N_{\lambda ,+})\sim\sum_{j'\in\Bbb N}\tilde c_{2j',+}(-\lambda
)^{\frac n2-j'-N},
$$
which contributes to the power $(-\lambda )^{-N}$ when $n$ is {\it
even}.
On the cylinder $X'\times \crp$, the Dirichlet s.g.o.\ $G_\lambda $
can be explicitly constructed (cf.\ (2.10)) to be$$
\gathered
G_\lambda =-K_{A_\lambda }(2A_\lambda )^{-1}T_{A_\lambda },\text{ where}\\
A_\lambda =(P'_1-\lambda )^{\frac12},\; K_{A_\lambda }\:
v(x')\mapsto e^{-x_nA_\lambda }v,\;T_{A_\lambda }\:u(x)\mapsto
\int_0^\infty e^{-x_nA_\lambda }u(x',x_n)\,dx_n;
\endgathered$$ with notation as in
\cite{G6}. Then $\tr_n G_\lambda =-(2A_\lambda
)^{-2}=-\frac14(P'_1-\lambda )^{-1}$. This is a resolvent on $X'$, and
the $N$'th derivative has a trace expansion
$$
 \Tr (G^{(N)}_\lambda )=\Tr_{X'}(\tr_n G_\lambda ^{(N)})
\sim\sum_{j'\in\Bbb N}a_{2j'}(-\lambda
)^{\frac {n-1}2-j'-N},
$$
which contributes to the power $(-\lambda )^{-N}$ when $n$ is {\it odd}.
\endexample

However, in  cases where $p$ has a parity that fits with the
dimension $n$, one can at least observe that the integral
$\int_{\Bbb R^n_+} \tslint \tr p(x,\xi )\,\d\xi dx$ has a coordinate
invariant meaning (cf\. Remark 3.2), so it gives an explicit value
that can be considered as the nonlocal part of
$C_0(P_+,P_{1,\operatorname{D}})$.
\medskip

Theorem 4.1 assures the validity of (2.7) for
$C_0(P_++G,P_{1,\operatorname{D}})$. In the proof of (2.8), we
restrict the attention to operators $P$ of normal order $\le 0$.

For (2.8),
it would be possible to refer again to Laguerre expansions (as in the
considerations of s.g.o.s in Section 3), using that $P$ on the symbol
level acts like a Toeplitz operator. More precisely, if one writes
$p(x',0,\xi )=\sum_{j\in\Bbb Z}a_j(x',\xi ')\hat\psi _j(\xi _n,\sigma
)$ with $$
\hat\psi _j(\xi _n,\sigma )=\frac{(\sigma -i\xi _n)^j}{(\sigma +i\xi
_n)^{j}},\tag4.7$$
then in the one-dimensional calculus based on Laguerre
expansion in $L_2(\rp)$, \linebreak$p(x',0,\xi ',D_n)_+$ acts like
the Toeplitz operator with matrix $(a_{j-k})_{j,k\in\Bbb N}$, cf\.
\cite{G2, Rem\. 2.2.12}. With some effort, the localness of
$C_0([P_+,G'],P_{1,\operatorname{D}})$ and
$C_0([P_+,P'_+],P_{1,\operatorname{D}})$ can be proved by use of the
localness of $C_0$ on commutators over $X'$, but the interpretations
and composition rules for the involved operators are not altogether
simple to deal with. Instead we shall rely on repeated commutator techniques.

In order to study the commutator of a pseudodifferential operator
 $P_+$ of order $\nu \in\Bbb Z$ with an
operator $A'=P'_++G'$ of order $\nu '\in \Bbb Z$, we write:$$
\Tr([P_+,A']R_\lambda ^N)=
\Tr([P_+,A']Q_{\lambda ,+} ^N)+\Tr([P_+,A']G_\lambda^{(N)}).\tag4.8
$$
We know from Proposition B.3 that the second term contributes only locally
to
$C_0$, and go on with an  analysis the first term:
$$\aligned
\Tr([P_+,A']&Q_{\lambda ,+} ^N)=\Tr(A'Q_{\lambda ,+}
^NP_+)-\Tr(A'P_+Q_{\lambda ,+} ^N)\\
&=\Tr(A'[Q_{\lambda } ^N,P]_+)-\Tr(A'G^+(Q_{\lambda } ^N)G^-(P))+
\Tr(A'G^+(P)G^-(Q_{\lambda } ^N))\\
&=\Tr(A'[Q_{\lambda } ^N,P]_+)-\Tr(G^-(P)A'G^+(Q_{\lambda } ^N))+
\Tr(A'G^+(P)G^-(Q_{\lambda } ^N)).
\endaligned\tag4.9$$
Again, the last two terms contribute only locally to $C_0$, by
Proposition B.3.
It is used that our hypotheses on $P$ assure that the s.g.o.s
$G^{\pm}(P)$ are of class 0.

Now consider the remaining term $\Tr(A'[Q_{\lambda } ^N,P]_+)$. The
expression $[Q_{\lambda } ^N,P]_+$ is
not very convenient in the boundary calculus, as
a mixture of strongly polyhomogeneous parameter-dependent
and arbitrary parameter-independent
interior operators.  However, we shall now show that it can be
reduced to an expression with
$\lambda $-dependent factors to the right only, with better
decrease in $\lambda $ than $Q_\lambda $, plus a manageable remainder
term.

\proclaim{Lemma 4.3} Let $P$ be of order $\nu $
and normal order $\le 0$.
For  $r>0$, let
$P_{(r)}$ denote the $r$'th
commutator of $P$ with $P_1$:
$$
P_{(1)}=[P,P_1],\;
P_{(2)}=[[P,P_1],P_1],\;\dots,\;P_{(r)}=[\cdots[[P,P_1],P_1]\dots,P_1];\tag4.10
$$
it is of order $\nu +r$ and of normal order
$\le \min\{\nu +r,2r\}$).

For any $M>0$, $[Q_\lambda ,P]$ may be written as$$
[Q_\lambda ,P]=P_{(1)}Q_\lambda ^2+P_{(2)}Q_\lambda
^3+\dots+P_{(M)}Q_\lambda ^{M+1}+Q_\lambda P_{(M+1)}Q_\lambda
^{M+1}.\tag4.11
$$
\endproclaim

\demo{Proof} Since $P_1$ has scalar principal symbol and is of order
2, $[P,P_1]$ is of order $\nu +1$. Since $P$ is of normal order 0 and
$P_1$ is of normal order 2,
$[P,P_1]$ is of normal order 2 (e.g., $\partial _{x_1}p_{1,2}$ is
generally so). Similarly, every subsequent
commutation lifts the order by 1 and the normal order by 2,
so the $r$'th commutator is of order $\nu +r$ and normal order $\le 2r$.
Since it has the transmission property, it also has normal order
$\le \nu +r$.
This explains the first statement. The second statement follows by
successive applications of the
following calculation:
$$\aligned
[Q_\lambda ,P_{(r)}]&=Q_\lambda P_{(r)}-P_{(r)}Q_\lambda
=Q_\lambda P_{(r)}(P_1-\lambda )Q_\lambda -Q_\lambda (P_1-\lambda
)P_{(r)}Q_\lambda \\
&=Q_\lambda [P_{(r)},P_1]Q_\lambda=Q_\lambda P_{(r+1)}Q_\lambda =
P_{(r+1)}Q_\lambda^2+ [Q_\lambda ,P_{(r+1)}]Q_\lambda;
\endaligned$$
hence
$$\aligned
[Q_\lambda ,P]&= P_{(1)}Q_\lambda ^2+[Q_\lambda ,P_{(1)}]Q_\lambda
=\dots \\
&=P_{(1)}Q_\lambda ^2+P_{(2)}Q_\lambda
^3+\dots+P_{(M)}Q_\lambda ^{M+1}+[Q_\lambda ,P_{(M)}]Q_\lambda ^{M}\\
&=P_{(1)}Q_\lambda ^2+P_{(2)}Q_\lambda
^3+\dots+P_{(M)}Q_\lambda ^{M+1}+Q_\lambda P_{(M+1)}Q_\lambda
^{M+1}.\quad\square
\endaligned$$
\enddemo

\proclaim{ Proposition 4.4} Let $P$ and $P'$ be $\psi $do's of order
$\nu $ resp\. $\nu '\in\Bbb Z$ and normal order $\le 0$, and let $G'$ be
an s.g.o.\ of order $\nu '$ and class $0$; denote $P'_++G'=A'$. For
$N>(\nu +\nu '+n)/2$,  $\Tr(A'[Q_\lambda ^N,P]_+)$
has an expansion$$
\Tr(A'[Q_\lambda ^N,P]_+)
\sim\sum_{j\ge
0}\!a_{j}(-\lambda)
^{\frac{n+\nu +\nu '- j}{2}-N}
+\sum_{k\ge
1}(a'_k\log (-\lambda )+a''_k)(-\lambda )
^{-\frac k2-N}.\tag 4.12
$$
\endproclaim

\demo{Proof} We know from the preceding calculation that there is an
expansion as in (1.10) with $\nu $ replaced by $\nu +\nu '$; the
point is to show that the series in $k$ starts with a lower power
than $-N$.

Recalling that $Q_\lambda ^N=\tfrac 1{(N-1)!}\partial _\lambda
^{N-1}Q_\lambda $, we find from (4.11):
$$
\multline
[Q^N_\lambda ,P]=c_1P_{(1)}Q_\lambda ^{N+1}+c_2P_{(2)}Q_\lambda
^{N+2}+\dots+c_MP_{(M)}Q_\lambda ^{N+M}\\
+\sum_{1\le l\le N}c'_{l}Q_\lambda ^{l}P_{(M+1)}Q_\lambda
^{N+M+1-l}.\endmultline\tag4.13
$$
In order to establish the expansion (4.12) we shall study the terms
$$P'_+(P_{(r)}Q^{(N+r)}_\lambda)_+ \
\text{ and }
\ G'(P_{(r)}Q^{(N+r)}_\lambda)_+; \quad 1\le r\le M; \tag4.14 $$
as well as
$$A'(Q_\lambda ^{l}P_{(M+1)}Q_\lambda ^{N+M+1-l})_+ ;\quad
1\le l\le N. \tag4.15$$
The operator $Q_\lambda ^{l}P_{(M+1)}Q_\lambda ^{N+M+1-l}$
is a $\psi$do on $\widetilde X$
with a symbol in $S^{\nu+M+1,0,-2N-2M-2}$.
The compositions in (4.15) will therefore be
trace-class on  $L^2(X)$ when $\nu -M-2N-1\le 0$, $\nu '+\nu -M-2N-1<-n$,
and their trace will
for any given $M'$ be $O(\lambda^{-M'})$ for sufficiently large $M$.
Our assertion will thus be true for these terms.

We turn to the terms in (4.14), considered in local coordinates.
In view of (2.6), we can decompose
$P_{(r)}= P_{(r)}'+P_{(r)}''$,
where $P_{(r)}''$ is a differential operator of order $\le \nu+r$ and
normal order $\le \min\{\nu +r,2r\}$, and
$P_{(r)}'$ has the property that $G^+(P'_{(r)})$ is of class 0, as a
sum of an operator of normal order 0 and an operator with a factor
$x_n^{M'}$ to the right, $M'\ge 2r$.
Then $$
(P_{(r)}Q_\lambda ^{N+r})_+
=P'_{(r),+}Q_{\lambda ,+}^{N+r}+G^+(P'_{(r)})G^-(Q_{\lambda
}^{N+r})+(P''_{(r)}Q_\lambda ^{N+r})_+. \tag4.16
$$

Note that
$$
P''_{(r)}=\sum_{0\le j\le 2r}S_j(x,D_{x'})D_{x_n}^j
$$
with tangential differential operators $S_j$ of order $\nu +r-j$.
The symbol of $Q^{N+r}_\lambda $ has the structure described
in (A.10), with typical term $r_{N+r,J,m}(p_{1,2}+\mu ^2)^{-m}$ where
$J+2N+2r\le 2m\le 4J+2N+2r$
 and $r_{N+r,J,m}$ is polynomial in $\xi $
of degree $\le 2m-2N-2r-J$. When
this is multiplied by $\xi _n^{2r}$ and we redefine $J+2r=J'$, we get
a term with numerator of order $\le 2m-2N-J'$ and denominator
$(p_{1,2}+\mu ^2)^{m}$, $J'+2N\le 2m\le 4J'+2N-6r\le 4J'+2N$; 
this is
of the type in the symbol expansion of $Q^N_\lambda $. Similarly,
multiplication by $\xi _n^{j}$ for $j\le 2r$ gives  terms of the types in
the
symbol expansions of $Q^{N'}_\lambda $ with $N'\ge N$.
The point of this analysis is that we can write
$$
G^-(P''_{(r)}Q^{N+r}_\lambda )=\sum_{0\le j\le
2r}S_j(x,D_{x'})G^-(D_{x_n}^jQ^{N+r}_\lambda ),
\tag4.17
$$
where the factors 
$G^-(D_{x_n}^{j}Q^{N+r}_\lambda )$
have structures like the 
$G^-(Q^{N'}_\lambda )$, $N'\ge N$, so that
Proposition B.3 can be applied.

The composition of $A'$ with the terms in (4.16) gives a number of terms:
$$\aligned
\text{(i) }&P'_+P'_{(r),+}Q_{\lambda ,+}^{N+r}=(P'P'_{(r)})_+Q_{\lambda
,+}^{N+r} -
G^+(P')G^-(P'_{(r)})Q_{\lambda ,+}^{N+r}\\
&=(P'P'_{(r)}Q_{\lambda
}^{N+r} )_+
-G^+(P'P'_{(r)})G^-(Q_{\lambda ,+}^{N+r})-G^+(P')G^-(P'_{(r)})Q_{\lambda
,+}^{N+r} ,\\
\text{(ii) }&G'P'_{(r),+}Q_{\lambda ,+}^{N+r},\\
\text{(iii) }&A'G^+(P'_{(r)})G^-(Q_{\lambda
}^{N+r}),\\
\text{(iv) }&P'_+(P''_{(r)}Q_\lambda ^{N+r})_+=
(P'P''_{(r)}Q_\lambda ^{N+r})_+-
G^+(P')G^-(P''_{(r)}Q_\lambda ^{N+r}),\\
\text{(v) }&G'(P''_{(r)}Q_\lambda ^{N+r})_+.
\endaligned$$
To the terms in (i) and (iii) with the s.g.o.-factor $G^-(Q_\lambda
^{r+N})$ we can apply
 Proposition B.3 with $\nu ,N$ replaced by $\nu +\nu ',N+r$; here the
series in $k$ starts with the power $(-\lambda )^{-\frac12-N-r}$.

The term $G^+(P')G^-(P''_{(r)}Q_\lambda ^{N+r})$ in (iv) with the
s.g.o.-factor $G^-(P''_{(r)}Q_\lambda ^{r+N})$ is written in view of
(4.17) as
$$
G^+(P')G^-(P''_{(r)}Q_\lambda ^{N+r})=
\sum_{0\le j\le
2r}G^+(P')S_j(x,D_{x'})G^-(D_{x_n}^jQ^{N+r}_\lambda) ,
$$
where each $G^+(P')S_j(x,D_{x'})$ is of
class 0. Proposition B.3 applies and gives contributions as in (4.12).

For the terms in (i) and (ii) of the form of an s.g.o.\ composed with
$Q^{r+N}_{\lambda ,+}$, we apply Proposition B.5. Since $r\ge 1$, the
series in $k$ has only powers $\le -1-N$ of $\lambda $, so the
contribution to the coefficient of $(-\lambda )^{-N}$ is local.

For (v),
we rewrite 
$G'(P''_{(r)}Q_\lambda ^{N+r})_+=(G'P''_{(r),+})Q_{\lambda ,+}
^{N+r}$,  
noting that the leftover term for the product is zero, 
since $P''_{(r)}$ is a differential operator. Now $G'P''_{(r),+}$
is an s.g.o.\ of order $\nu'+\nu+r$ and class $\le\min\{\nu+r, 2r\}$. It
therefore has a representation 
$$G_0 + \sum_{0\le j\le 2r-1} K_j\gamma_j$$ 
with an s.g.o.\ $G_0$ of order $\nu+\nu'+r$ and class 0 and 
Poisson operators $K_j$ of order $\nu+\nu'+r-j$.
The composition $G_0Q_{\lambda ,+}^{N+r}$
is like (ii); Proposition B.5 
applies to give an expansion (4.12).

Next we decompose $K_j$ according to Lemma B.7:
$K_j=\sum_{l\ge 0}\Phi_lC_{jl}$
with a rapidly decreasing sequence
$(C_{jl})_l$  in $S^{\nu+\nu'+r-j+\frac12}$. Then
$$\Tr(K_j\gamma_jQ^{N+r}_{\lambda,+})=
\sum_{l\ge 0}
\Tr_{\Bbb R^{n-1}}(\gamma_0\partial^j_{x_n}
Q^{N+r}_{\lambda ,+}\Phi_lC_{jl}).$$ 
The composition $\gamma_0\partial_{x_n}^jQ^{N+r}_{\lambda ,+}\Phi_l$
is a $\psi$do on $\Bbb R^{n-1}$
with the symbol
$$\sum_{j_1+j_2=j}\binom{j}{j_1}
\int^+
(i\xi_n)^{j_1}\,\partial^{j_2}_{x_n}q^{\circ(N+r)}(x',0,\xi,\mu)\,
\hat\varphi_l(\xi_n,\sigma)\d\xi_n.\tag 4.18$$
We know from (A.8)--(A.11) that 
$(i\xi_n)^{j_1}\partial_{x_n}^{j_2}q^{\circ(N+r)}_{-2(N+r)-J}$ is
a finite sum of terms of the form $\xi_n^{j_1}r_{J,m}(p_{12}+\mu^2)^{-m}$
with a polynomial $r_{J,m}$ of degree $\le 2m-2(N+r)-J$. 
De\-com\-po\-sing it into simple 
fractions, we obtain a sum of terms of the form 
$r^\pm_{J,k}(\kappa^\pm\pm i\xi_n)^{-k}$ with
$r^\pm_{J,k}(x',\xi',\mu)$ strongly 
homogeneous in $(\xi',\mu)$ of degree $\le j+k-2(N+r)-J$. 
Inserting this into (4.18) and using Lemma B.1, 
we see that 
$\gamma_jQ^{N+r}_{\lambda ,+}\Phi_l$  has its symbol in
$S^{\frac12,0,-2N-2r+j}$; the symbol seminorms grow at most 
polynomially in $l$. 
Since $j<2r$, we then reach the expansion (4.12).

The remaining $\psi $do terms in (i) and (iv), $(P'P'_{(r)}Q_{\lambda
}^{N+r} )_+$ and $(P'P''_{(r)}Q_\lambda ^{N+r})_+$, add up to
$(P'P_{(r)}Q_\lambda ^{N+r})_+$, which is treated as in Theorem 4.1,
using the pointwise kernel expansion on $\widetilde X$.
Here $P'P_{(r)}Q_\lambda ^{r+N}$ has symbol in $S^{\nu '+\nu +r,
0,-2r-2N}\subset S^{\nu '+\nu -r-2N,0,0}$\linebreak$\cap S^{\nu '+\nu +r,
-2r-2N,0}$ with $d$-index $\le -2-2N$ since $r\ge 1$, so we get an
expansion
contributing only locally to the coefficient of $(-\lambda )^{-N}$.
\qed
\enddemo

We now have all the ingredients for the proof of:

\proclaim{Theorem 4.5}
Let $P$ and $P'$ be $\psi $do's of order
$\nu $ resp\. $\nu '\in\Bbb Z$ and normal order $0$, and let $G$ and $G'$
be
singular Green operators of order $\nu $ resp.\ $\nu '$ and class
$0$. Then $C_0([P_++G,P'_++G'],P_{1,\operatorname{D}})$ is locally
determined.
\endproclaim

\demo{Proof} We have that
$$[P_++G,P'_++G']=[P_++G,P'_+]+[P_+,G']+[G,G'].
$$
The last term was shown in Theorem 3.13 to contribute locally to
$C_0$.
For the two other terms we
have the analysis above, through (4.8), (4.9) and Proposition 4.4,
showing that they contribute only locally to $C_0$.
\qed
\enddemo

This completes the proof that (2.8) holds also for
operators of type $A=P_++G$,
so
$C_0(P_++G,P_{1,\operatorname{D}})$ is
indeed a quasi-trace on such operators. In addition, we have found
the interesting information that the contribution from $P_+$ can be
traced back to  a pointwise defined contribution from $P$ over
$\widetilde X$, and that the contribution from $G$ can
be traced back to an interior trace contribution of order $-\infty $
plus a contribution from the
normal trace on $X'$; here both $\widetilde X$ and $X'$ are
compact manifolds without boundary.

\head Appendix A. The structure of the auxiliary operators
\endhead

We here recall the symbol formulas established and used in \cite{GSc},
and some useful consequences.
 The parameter-dependent entries were
indexed by $\mu =(-\lambda )^{\frac12}$ in \cite{GSc}; we simply replace
this here by
indexation by $\lambda $, although $\mu $ can still appear as
a variable.

Throughout this paper,
we denote by $[\xi ']$ a positive $C^\infty $ function of $\xi
'$ that coincides with $|\xi '|$ for $|\xi '|\ge 1$. It
will often be denoted $\sigma (\xi ')$ or just
$\sigma $.

The principal (second-order) symbol of $P_1$ is denoted $p_{1,2}$, so
the principal symbol of
$Q_\lambda $ is $q_{-2}=(p_{1,2}-\lambda)^{-1}
=(p_{1,2}+\mu ^2)^{-1}$.

We assume $P_1$
to be strongly elliptic; this means that the principal symbol
$p_{1,2}(x,\xi )$ has positive real part when $\xi \ne 0$. Since
$|\operatorname{Im}p_{1,2}(x,\xi )|\le |p_{1,2}(x,\xi )|\le
C\operatorname{Re}p_{1,2}(x,\xi )$, there
is a sector $\Gamma $ such that
$p_{1,2}(x,\xi )+\mu ^2\ne 0$ when  $\mu \in \Gamma \cup \{0\}$, $(\xi
,\mu
)\ne (0,0)$; cf\. also (3.9).

The following observation will be useful:

\proclaim{Lemma A.1} Let $p (x,\xi )$ be a uniformly strongly
elliptic homogeneous second-order differential operator symbol
on $\Bbb R^n$,
with $p-\lambda =p+\mu ^2$ invertible for $\mu \in\Gamma $.
Then for any $m$, $L>0$,$$
\aligned
(p +\mu ^2)^{-m}&=
\sum_{0\le j <L}c_{m,j } \mu ^{-2m-2j } p^{j }+
\mu ^{-2m-2L } p'_L(x,\xi ,\mu ), \text{ where}\\
 p'_L(x,\xi ,\mu )&=\sum_{0\le k\le m}c_{m,L,k}
 p ^{L+k}(p +\mu ^2)^{-k}. \endaligned\tag A.1$$
Here one has for all indices:$$
|\partial _\lambda ^N\partial _x^\beta \partial _\xi ^\alpha
p'_L(x,\xi ,\mu )|\le C_{\alpha ,\beta ,N
}\ang\xi ^{2L-|\alpha |-2N},\tag A.2
$$
with constants $C_{\alpha ,\beta ,N}$
independent of $\mu \in\Gamma $ for $|\mu |\ge 1$.
\endproclaim

\demo{Proof}
One has for $b\in \Bbb C\setminus\{1\}$, $L\ge 1$,
$$\aligned (1-b)^{-m}=&\frac1{(m-1)!}\partial_b^{m-1}\frac1{1-b}=
\frac1{(m-1)!}\partial_b^{m-1}\Bigl(\sum_{0\le k <L+m-1}b^k
+\frac{b^{L+m-1}}{1-b}\Bigr)
\\
=&\sum_{0\le j <L}c_j b^{j }+\sum_{1\le l\le m}
\frac{c'_lb^{L+l-1}}{(1-b)^{l}}.
\endaligned\tag A.3$$
An application to $1-b=1+p /\mu ^2$ gives$$
\aligned
(p +\mu ^2)^{-m}&=\mu ^{-2m}(1+p /\mu ^2)^{-m}\\
&=\mu ^{-2m}\bigl(\sum_{0\le j <L}c_j \mu ^{-2j
}p ^j +\sum_{1\le l\le m}c'_l\mu
^{-2L-2l+2}p ^{L+l-1}
(1+p /\mu ^2)^{-l}\bigr)\\
&=\sum_{0\le j <L}c_j  \mu ^{-2m-2j
}p ^j +\mu ^{-2m-2L+2}\sum_{1\le l\le m}c'_lp ^{L+l-1}
(p +\mu ^2)^{-l}.
\endaligned$$
Insertion of $\mu ^2(p +\mu ^2)^{-1}=1-p (p +\mu ^2)^{-1}$
in the last sum leads to (A.1). Since
$p $
is polynomial of degree 2,
$(p +\mu ^2)^{-1}$ is $C^\infty $ and homogeneous
in $(\xi ,\mu
)$ of degree $-2$ for $(\xi ,\mu )\ne 0$, and $\partial
_\lambda ^N(p-\lambda )^{-l}=c_{l,N}(p-\lambda )^{-l-N},$
$$
\multline
|\partial _\lambda ^N\partial _x^\beta
\partial _\xi ^\alpha p'_L(x,\xi ,\mu )|
\le C'_{\alpha ,\beta ,N}\sum_{0\le k\le m,\alpha '\le\alpha }
(\partial _\xi ^{\alpha '}p^{L+k})
|(\xi ,\mu )|^{-2k-2N-|\alpha -\alpha '|}\\
\le C_{\alpha ,\beta,N}
\ang\xi ^{2L-|\alpha |-2N}\endmultline
$$
for $|\mu |\ge 1$, showing (A.2).
\qed
\enddemo

Let $\xi =(\xi ',\xi _n)$ in local coordinates
at  the boundary, then for $(\xi ',\mu )\ne (0,0)$,
the strong ellipticity implies that the polynomial in
$\xi _n$
$$
p_{1,2}(x',0,\xi ',\xi _n)+\mu ^2=a(x')\xi _n^2+b(x',\xi ')\xi
_n+c(x',\xi ')+\mu ^2\tag A.4
$$ has two roots $\varrho _1(x',\xi ',\mu )$ and
$\varrho _2(x',\xi ',\mu )$ in $\Bbb C\setminus \Bbb R$.
When $\mu \in\crp$, one
of the roots, say $\varrho _1$, lies in $\Bbb C_+$ and the other, $\varrho
_2$, in $\Bbb C_-$. For, $a\xi _n^2+b\xi _n+c+\mu ^2$ can be carried
into  $\operatorname{Re}a\,\xi _n^2+\operatorname{Re}b\,\xi
_n+\operatorname{Re}c+\mu ^2$ by a homotopy that preserves the property
of having positive real part, and for the latter polynomial, the roots are
placed in this way; they depend continuously on the
polynomial, hence cannot cross the real axis. The placement of the
roots is also
preserved when $\mu $ is moved to a general element of $\Gamma $.
Thus we can denote the roots $\pm i\kappa ^{\pm}(x',\xi ',\mu )$,
 where
$\kappa ^\pm$ have positive real part; they depend smoothly on
$(x',\xi ',\mu )$ for $(\xi ',\mu )\ne (0,0)$. They are homogeneous of
degree
1 in $(\xi ',\mu )$ and bounded away from $\Bbb R$ for $|(\xi ',\mu
)|=1$, $\mu \in \Gamma \cup \{0\}$, so in fact they take values in
a sector $\{|\arg z|\le \frac \pi 2-\delta \}$ for some $\delta >0$.

We recall from \cite{GSc, Sect\. 2.b} (with some small precisions):

\proclaim{Lemma A.2} The symbol of $Q_\lambda $ has the following form in
local coordinates, for $\mu \in \Gamma $:
$$\aligned
q(x,\xi ,\mu )&\sim \sum_{l\in\Bbb N}q_{-2-l}(x,\xi ,\mu ),\text{ with}\\
q_{-2}(x,\xi ,\mu )&=(p_{1,2}(x,\xi )+\mu ^2)^{-1}I,\\
q_{-2-J}(x,\xi ,\mu )&=\sum_{J/2+1\le m\le 2J+1}\frac{r_{J,m}(x,\xi
)}{(p_{1,2}(x,\xi )+\mu ^2)^{m}},\text{ for }J\ge 0;\endaligned\tag A.5
$$
here the $r_{J,m}$ are $(\operatorname {dim} E\times \operatorname{dim} E)$-matrices of
homogeneous polynomials in $\xi $ of degree $2m-2-J$ with smooth
coefficients, and the remainders $q^{\prime}_{-2-M}=q-\sum_{0\le
J<M}q_{-2-J}$
satisfy estimates for all indices $\partial _{x}^\beta \partial _\xi
^\alpha
\partial _\mu ^kq^{\prime}_{-2-M}=O(\ang{\xi ',\mu
}^{-2-M-|\alpha |-k})$, for $|\mu |\ge 1$, $\mu $ in closed
subsectors of $\Gamma $.

Concerning $q(x',0,\xi )$, we have: Writing
$$\aligned
p_{1,2}(x',0,\xi )+\mu ^2&=a(x')(\xi _n-i\kappa ^+(x',\xi ',\mu ))
(\xi _n+i\kappa ^-(x',\xi ',\mu ))\\
&=
a(x')(\kappa ^+(x',\xi ',\mu )+i\xi _n)
(\kappa ^-(x',\xi ',\mu )-i\xi _n),\endaligned\tag A.6
$$
we can decompose each term
in simple fractions (at $x_n=0$):
$$\aligned
\frac{r_{J,m}(x',0,\xi )}{(p_{1,2}(x',\xi )+\mu
^2)^{m}}&=h^+\frac{r_{J,m}}{(p_{1,2}+\mu
^2)^{m}}+h^-\frac{r_{J,m}}{(p_{1,2}+\mu ^2)^{m}},\\
h^\pm \frac{r_{J,m}}{(p_{1,2}+\mu ^2)^{m}}&=
\sum_{1\le j\le m}\frac{r^\pm_{J,m,j}(x',\xi ',\mu
)}{(\kappa ^\pm(x',\xi ',\mu )\pm i\xi
_n)^{j}},\endaligned\tag A.7
$$
where the $r_{J,m,j}^\pm(x',\xi ',\mu )$ are strongly homogeneous of
degree
$j-J-2$ in $(\xi ',\mu )$.
This gives a decomposition of
$q(x',0,\xi ',\mu )$:
$$
q=h^+q+h^-q=q^{+}+q^{-},\quad
q^{\pm}(x',\xi ,\mu )\sim\sum_{
J\ge 0}q^{ \pm}_{-2-J}(x',\xi ,\mu ),
$$with terms as in {\rm (A.7)}.

The structure of the normal derivatives is similar:
We have
$$
\partial _{x_n}^lq_{-2-J}(x,\xi ,\mu )
=\sum_{1+J/2\le m\le 2J+1+l}
\frac{r^{l}_{J,m}(x,\xi)}{(p_{1,2}(x,\xi )+\mu ^2)^{m}},
$$
for all $l$, with homogeneous polynomials $r^{l}_{J,m}(x,\xi )$
of degree $2m-2-J$ in $\xi $. At $x_n=0$ we can decompose as before,
obtaining
$$\aligned
\quad \partial _{x_n}^lq=\partial _{x_n}^lq^++\partial _{x_n}^lq^{-},&
\quad
\partial _{x_n}^lq^{\pm}(x',\xi ,\mu )
\sim\sum_{J\ge 0}
\partial _{x_n}^lq^{\pm}_{-2-J}(x',\xi ,\mu ),\\
\partial _{x_n}^lq_{-2-J}(x',\xi ,\mu )
&=\partial _{x_n}^lq^{+}_{-2-J}(x',\xi ,\mu
)+\partial _{x_n}^lq^{-}_{-2-J}(x',\xi ,\mu ),\\
\partial _{x_n}^lq^{\pm}_{-2-J}(x',\xi ,\mu )
&=\sum_{1\le j\le 2J+1+l}
\frac{r^{l,\pm}_{J,j}(x',\xi ',\mu )}
{(\kappa ^\pm(x',\xi ',\mu )\pm i\xi_n)^{j}},
\endaligned\tag A.8
$$
where
the numerators $r^{l,\pm}_{J,j}(x',\xi ',\mu )$ are strongly
homogeneous of degree $j-J-2$.
\endproclaim

It is useful to observe that the $r^{\pm}_{J,m,j}$
as well as $\kappa ^{\pm}$ are in fact functions of $(x',\xi
',\lambda )$, $\lambda =-\mu ^2$, strongly quasihomogeneous in $(\xi
',\lambda )$ with weight $(1,2)$ in the following sense: We say that
$r(x',\xi ',\lambda )$ is $(1,2)$-homogeneous in $(\xi ',\lambda )$
of degree $d$, when
$$
r(x',t\xi ',t^2\lambda )=t^{d}r(x',\xi ',\lambda );\tag A.9
$$
it is strongly so if (A.9) holds for $|\xi '|+|\lambda |\ge \varepsilon $,
weakly so if it holds for $|\xi '|\ge \varepsilon $. Then since
$r(x',\xi ',\lambda )=t^{-d}r(x',t\xi ',t^2\lambda )$, it is
readily checked that $\partial _\lambda ^Nr$ is
$(1,2)$-ho\-mo\-ge\-ne\-ous of
degree $d-2N$. Thus for the symbols that depend on $\mu $ through
$\lambda $ in this way, differentiation with respect to $\lambda $
lowers the homogeneity degree in $(\xi ',\mu )$ by two steps,
preserving strong homogeneity.
The estimates of $\lambda $-derivatives of the remainders
$q^{\prime}_{-2-M}$
likewise improve by two steps for each derivative.

These considerations of $\lambda $-derivatives, playing on the
strong quasi-ho\-mo\-ge\-ne\-ity of symbols coming from
$R_\lambda $, will replace the calculations for
higher powers of $R_\lambda $ used in \cite{GSc}, cf\. (2.12). (The
calculus
in \cite{GH} is set up to handle the anisotropic homogeneity in
terms of $\xi '$ and $\lambda $ and could give further information on
higher terms in the full trace expansions, but for the discussion of
the leading
nonlocal term we carry out here, the calculus of \cite{GS1}
will suffice.) For convenience, we recall explicitly the structure of
the formulas for
symbols connected with higher powers of $Q_\lambda $, denoting the
symbol of
$Q_\lambda ^N=\frac{\partial _\lambda^{N-1}}{(N-1)!}Q_\lambda $
by $q^{\circ N}$:
$$\aligned
q^{\circ N}(x,\xi ,\mu )
&\sim \sum_{l\in\Bbb N}q^{\circ N}_{-2N-l}(x,\xi ,\mu ),
\text{ with }q^{\circ N}_{-2N}=(p_{1,2}+\mu ^2)^{-N},\\
q^{\circ N}_{-2N-J}(x,\xi ,\mu )&=\sum_{J/2+N\le m\le 2J+N}
\frac{r_{N,J,m}(x,\xi)}{(p_{1,2}(x,\xi )+\mu ^2)^{m}},
\text{ for }J\ge 0,
\endaligned\tag A.10
$$
where the $r_{N,J,m}$ are
homogeneous polynomials in $\xi $ of degree $2m-2N-J$
with smooth coefficients. Moreover for $x=(x',0)$,
$$
\aligned
q^{\circ N}&=h^+q^{\circ N}+h^-q^{\circ N}=q^{\circ N,+}+q^{\circ N,-},\\
q^{\circ N,\pm}(x',\xi ,\mu )&\sim\sum_{
J\ge 0}q^{\circ N, \pm}_{-2-J}(x',\xi ,\mu ),\\
q^{\circ N, \pm}_{-2-J}(x',\xi ,\mu )&=\sum_{1\le j\le
2J+N}\frac{r^\pm_{N,J,j}(x',\xi ',\mu
)}{(\kappa ^\pm(x',\xi ',\mu )\pm i\xi
_n)^{j}};
\endaligned\tag
A.11
$$
here $r_{N,J,j}^\pm(x',\xi ',\mu )$ are strongly homogeneous of degree
$j-J-2N$ in $(\xi ',\mu )$.

For 
$\partial_{x_n}^l q^{\circ N}(x,\xi ,\mu )$ and
$\partial_{x_n}^l q^{\circ N,\pm}(x',\xi ,\mu )$, we obtain
corresponding results: Their structure is as in 
(A.10) and (A.11), except
that now the summation will be over the sets $J/2+N\le m\le 2J+N+l$ and
$1\le j\le 2J+N+l$.
\medskip

The following symbols derived from the principal symbol of $q$ at
$x_n=0$ played an important role in \cite{GSc, Sect\. 5}:
$$\aligned
\alpha ^{(1)}(x',\xi ',\mu )&=[h^+q_{-2}(x',0,\xi ,\mu )]_{\xi _n=-i\sigma
}+
[h^-q_{-2}(x',0,\xi ,\mu )]_{\xi _n=i\sigma }\\
&=q^+_{-2}(x',\xi ',-i\sigma ,\mu )+q^-_{-2}(x',\xi ',i\sigma ,\mu
)\\
&=\frac 1{(\kappa ^++\kappa ^-)(\kappa ^++\sigma )}+\frac 1{(\kappa
^++\kappa ^-)(\kappa ^-+\sigma )},\\
\alpha ^{(N)}(x',\xi ',\mu )&=[h^+q^{\circ N}_{-2N}(x',0,\xi ,\mu )]_{\xi
_n=-i\sigma }+
[h^-q^{\circ N}_{-2N}(x',0,\xi ,\mu )]_{\xi _n=i\sigma }\\
&=\tfrac{\partial _\lambda ^{N-1}}{(N-1)!}\alpha ^{(1)}(x',\xi ',\mu );
\endaligned\tag A.12$$
here
$\sigma =[\xi ']$, and $\alpha ^{(N)}$ is weakly polyhomogeneous
in $S^{-2N,0,0}\cap S^{0,-2N,0}$ (in fact in $S^{0,0,-2N}$, see
below). The crucial
information established in \cite{GSc, Lemma 5.5} was that the
coefficient of $\mu
^{-2N}$ in the expansion in powers of $\mu $ (as in (3.4)) is 1:
$$
\alpha ^{(N)}(x',\xi ',\mu )=\mu ^{-2N}+
\alpha ^{(N)}_1(x',\xi ',\mu )\text{ with } \alpha ^{(N)}_1(x',\xi ',\mu
)\in S^{1,-2N-1,0};
\tag A.13
$$
the remainder $\alpha _1^{(N)}$ is also in $S^{0,-2N,0}$ since
$\alpha ^{(N)}$ and $\mu ^{-2N}$ are so.

\proclaim{Lemma A.3} The symbols $(\kappa ^\pm+\sigma )^{-1}$ are
weakly polyhomogeneous, belong to $S^{0,0,-1}$ and have $N$'th
$\lambda $-derivatives in $S^{0,0,-1-2N}$. For each $N$, $\alpha ^{(N)}$
is
weakly polyhomogeneous lying in $S^{0,0,-2N}$.
\endproclaim

\demo{Proof} First note that
$\kappa ^\pm$ and $(\kappa ^\pm)^{-1}$ are strongly polyhomogeneous in
$(\xi ',\mu )$ of degree $1$ resp\. $-1$, so they lie in
$S^{0,0,1}$ resp\. $S^{0,0,-1}$. As noted above, they are strongly
$(1,2)$-homogeneous in $(\xi ',\lambda )$ of degree 1 resp\. $-1$;
hence the $N$'th $\lambda $-derivatives lie in $S^{0,0,1-2N}$ resp\.
$S^{0,0,-1-2N}$. Similarly, since $\kappa ^+$ and $\kappa ^-$ lie in
a proper subsector of $\{z\in\Bbb C\mid \operatorname{Re}z>0\}$,
$|\kappa ^++\kappa ^-|\ge \operatorname{Re}(\kappa ^++\kappa ^-)\ge
c_0|(\xi ',\mu )|$ with $c_0>0$, so also
$(\kappa ^++\kappa ^-)^{-1}$ belongs to
$S^{0,0,-1}$, with $N$'th $\lambda $-derivatives in $S^{0,0,-1-2N}$.

It was observed in \cite{GSc} that $(\kappa ^\pm+\sigma )^{-1}\in
S^{-1,0}\cap
S^{0,-1}$; we now recall the proof showing how it also assures
that $(\kappa ^\pm+\sigma )^{-1}\in S^{0,0,-1}$:

Let $\kappa $ stand for $\kappa ^+$ or $\kappa ^-$.
Write $(\kappa  +\sigma
)^{-1}=(\kappa  )^{-1}(1+\sigma /\kappa  )^{-1}$. We know already
that $(\kappa
 )^{-1}\in S^{0,0,-1}$, so it remains to show that $(1+\sigma /\kappa
 )^{-1}$ is in $S^{0,0,0}$. Since $\sigma \in S^{1,0,0}$ and
$(\kappa  )^{-1}\in
S^{0,0,-1}\subset S^{-1,0,0}$,
$\sigma /\kappa  $ is in $S^{0,0,0}$, and so is $1+\sigma
/\kappa  $. Moreover, it is
bounded in $(\xi ',\mu )$ (for $\mu \in\Gamma $, $|\mu |\ge 1$).
Since $\operatorname{Re}\kappa  >0$,
$\sigma <\sigma +\operatorname{Re}\kappa  \le |\sigma +\kappa  |$,
so also the inverse $(1+\sigma
/\kappa  )^{-1}=\kappa  /(\kappa  +\sigma )=1-\sigma /(\sigma
+\kappa  )$ is bounded. Then \cite{GS1, Th\.
1.23} (just the beginning of the proof) shows that the inverse $(1+\sigma
/\kappa  )^{-1}$ does belong to $S^{0,0,0}$.

For
the $\lambda $-derivatives, we observe that
$\partial _\lambda (\kappa  +\sigma )^{-1}=
- (\kappa   +\sigma )^{-2}\partial _\lambda \kappa  \in
S^{0,0,-3}$ by the composition rules, and hence by successive application,
$\partial _\lambda ^k(\kappa  +\sigma )^{-1}\in S^{0,0,-1-2k}$ for
all $k$.

Using these informations, it is now seen from the form of $\alpha
^{(1)}$ in (A.12) that it lies
in $S^{0,0,-2}$ with $k$'th $\lambda $-derivatives in
$S^{0,0,-2-2k}$, so $\alpha ^{(N)}\in S^{0,0,-2N}$.
\qed
\enddemo

Besides $G^{(N)}_\lambda $ (cf\. (2.11)), we also need to consider
the s.g.o.s $G^{\pm}(Q^N)$, which arise from compositions such as
$$
P_+Q^N_{\lambda ,+}=(PQ^N_\lambda )_+-G^+(P)G^-(Q^N_\lambda ).\tag A.14
$$
Their symbols have the following structure:

\proclaim{Lemma A.4}
The operators $G^{(N)}_\lambda $ and $G^{\pm}(Q^N)$
have
symbols of the form$$\aligned
g(x',\xi ,\eta _n,\mu )&\sim
\sum_{J\in\Bbb N}g_{-2N-1-J}(x',\xi ,\eta _n,\mu ),\text{ with }\\
g_{-2N-1-J}(x',\xi ,\eta _n,\mu )&=\sum_{ \Sb j\ge 1,\,j'\ge 1\\
j+j'\le 2J+N+1\endSb}\tfrac{s_{J,j,j',N}(x',\xi ',\mu )}
{(\kappa +i\xi
_n)^{j}(\kappa '-i\eta
_n)^{j'}};
\endaligned\tag A.15$$
here $(\kappa ,\kappa ')$ equals $(\kappa ^+,\kappa ^-)$ for
$G^{(N)}_\lambda $, $(\kappa ^+,\kappa ^+)$ for $G^{+}(Q^N)$, and
$(\kappa ^-,\kappa ^-)$ for $G^{-}(Q^N)$.
The numerators $s_{J,j,j',N}(x',\xi ',\mu )$ are
strong\-ly homogeneous in $(\xi ',\mu )$ of degree
$j+j'-J-2N-1$.
They are in fact functions of $(x',\xi
',\lambda )$, $\lambda =-\mu ^2$, strongly $(1,2)$-homogeneous in $(\xi
',\lambda )$ of the indicated degrees, such that differentiation
with respect to $\lambda $ gives a strongly $(1,2)$-homogeneous
symbol of $2$ steps lower degree.

For the remainders $g'_{-2N-1-M}=g-\sum_{0\le
J<M}g_{-2N-1-J}$,  the sup-norms in $\xi _n$ resp\. $(\xi
_n,\eta _n)$ are bounded by $\ang{\xi ',\mu }$ in powers
$-2N-1-M$, decreasing by $|\alpha |$ for
differentiations in $(\xi ,\eta _n)$ of order $\alpha $,
and by $2$ for each differentiation in $\lambda $ (no
change of the power for differentiations in $x'$).

\endproclaim

\demo{Proof} The formulas for the symbols of $G^{(N)}_\lambda $ and
$G^\pm(Q^N_\lambda )$ come from
\cite{GSc, Prop\. 2.5 and 2.3}.
The remainder estimates hold since the symbols are strongly
polyhomogeneous in $(\xi ,\eta _n,\mu )$ so that
standard estimates hold when $|\mu |$ is considered as an extra
cotangent variable (on each ray). The
strong (1,2)-polyhomogeneity $(\xi
,\eta _n,\lambda )$ assures the statements on $\lambda$-derivatives. \qed
\enddemo

\head Appendix B. Formulas for compositions with Laguerre
functions
\endhead

Recall the formulas for the (Fourier transformed) Laguerre functions
we use in expansions
of parameter-independent operators:
$$
\hat\varphi '_k(\xi _n,\sigma )=\frac{(\sigma -i\xi _n)^k}{(\sigma
+i\xi _n)^{k+1}},\quad\hat\varphi _k(\xi _n,\sigma )
=(2\sigma )^{\frac12}\frac{(\sigma -i\xi _n)^k}{(\sigma
+i\xi _n)^{k+1}};\tag B.1
$$
A priori, $\sigma$ here might be any positive number, but we will take
$\sigma=[\xi']$.
The $\hat\varphi _k$ are the {\it normalized} Laguerre
functions, which for $k\in\Bbb Z$ form an orthonormal basis of
$L_2(\Bbb R)$.
Differentiations in $\xi '$ and $\xi _n$ follow the rules$$
\aligned
\partial_{\xi_j}\hat\varphi_k(\xi_n,\sigma)
&=\big(k\hat\varphi_{k-1}-
\hat\varphi_{k}
-(k+1)
\hat\varphi_{k+1}\big)(2\sigma)^{-1}\partial_{\xi_j}\sigma ,\quad
j<n,\\
\partial_{\xi_n}\hat\varphi_k(\xi_n,\sigma)
&=-i\big(k\hat\varphi_{k-1}+
(2k+1)\hat\varphi_{k}
+(k+1)
\hat\varphi_{k+1}\big)(2\sigma)^{-1} ;
\endaligned
\tag B.2
$$
here $\partial _{\xi _j}\sigma =\xi _j\sigma ^{-1}$ for $|\xi '|\ge 1$.

\medskip
We have from \cite{GSc} the following formulas for $\xi _n$-compositions
of Laguerre funtions and rational
functions involving $\kappa ^\pm$:

\proclaim{Lemma B.1} For all $m\ge 0$ and $j\ge 1$:
$$
\int \frac{(\sigma \pm i\xi _n)^m}{(\sigma \mp i\xi _n)^{m+1}}\frac{1}
{(\kappa ^\pm(x',\xi ',\mu )\pm i\xi
_n)^{j}}\,\d\xi _n=\sum_{m'\ge 0,\, |m'-m|<j}(\kappa^\pm) ^{1-j}
a^\pm_{jm'}\frac{(\sigma -\kappa ^\pm)^{m'}}{(\sigma
+\kappa ^\pm
)^{m'+1}},
\tag B.3$$
with universal constants $a^\pm_{jm'}$ that are $O(m^j)$ for fixed $j$.
The
resulting expressions are
weakly polyhomogeneous $\psi $do symbols belonging to
$S^{0,0,-j}(\Bbb R^{n-1}\times\Bbb R^{n-1},\Gamma )$, with $\partial
_\lambda ^N$ of the symbols lying in $S^{0,0,-j-2N}(\Bbb
R^{n-1}\times\Bbb R^{n-1},\Gamma )$.
\endproclaim

\demo{Proof}
The formulas are shown in \cite{GSc, Lemma 3.2}. The statements on
symbol classes follow since $\sigma =[\xi ']\in S^1\subset
S^{1,0,0}$, $\kappa ^{\pm}\in S^{0,0,1}$ with $N$'th $\lambda
$-derivatives in $S^{0,0,1-2N}$, and $(\sigma
+\kappa ^{\pm})^{-1}\in S^{0,0,-1}$ with $N$'th $\lambda
$-derivatives in $S^{0,0,-1-2N}$, by Lemma A.3. For the derivatives of the
composed expressions one can use the formulas
\cite{GSc, (3.21)} for the derivatives of $(\sigma -\kappa
^\pm)^m(\sigma +\kappa ^\pm)^{-m}$.
\qed\enddemo

These formulas enter in calculations of compositions such as
$GG^{(N)}_\lambda $ and $GG^{\pm}(Q^N_\lambda )$, where the symbol
of $G$ is Laguerre
expanded and the rational structure of the symbols of $G^{(N)}_\lambda
$ and $G^\pm(Q^N_\lambda )$ (A.15) is used.

\proclaim{Lemma B.2} The symbol
$$
s_{j,j',m}=
\int \frac{(\sigma -i\xi _n)^m}{(\sigma
+i\xi _n)^{m+1}}\frac
{1}{(\kappa ^++ i\xi
_n)^{j}(\kappa ^--
i\xi
_n)^{j'}}\,\d\xi _n
$$
satisfies for $m\ge 0$, $j$ and $j'\ge 1$:
$$s_{j,j',m}=\sum_{\Sb |m-m'|\le j''< j'\\  m'\ge
0\endSb}b_{jj'j''m'}(\kappa ^-)^{-j''}\frac{(\sigma -\kappa
^-)^{m'}}{(\sigma +\kappa ^-)^{m'+1}}(\kappa ^++\kappa ^-)^{-j-j'+1+j''},
\tag B.4$$
where the $b_{jj'j''m'}$ are universal constants that are
$O(m^{j'})$ for fixed $j,j'$.
This is a weakly polyhomogeneous symbol
in $S^{0,0,-j-j'}$, with $N$'th $\lambda $-derivatives in
$S^{0,0,-j-j'-2N}$, for all $N$.
There is a similar formula for $m\le 0$, with $j$ and $j'$
interchanged, $\kappa ^+$ and $\kappa ^-$
interchanged. \endproclaim

\demo{Proof}
The formulas were shown in \cite{GSc, Lemma 4.2}, and the symbols are
analyzed as in the preceding proof.
\qed\enddemo

These formulas enter in calculations of compositions such as
$P_+G^{(N)}_\lambda $ and $P_+G^\pm(Q^N_\lambda )$, where
the symbol of $P$ is Laguerre expanded and the rational structure of
the s.g.o.\ is used.

Composition formulas involving $\psi $do's on the interior
contain $h^+_{\xi _n}$- and $h^-_{\xi _n}$-pro\-ject\-ions. We recall
from \cite{GSc, (3.9), (4.8)} that the projections can be removed
in certain integrals, e.g.: When $\tilde q$ is a rational function of $\xi
_n$ of the form $r(x',\xi ',\mu )(\kappa ^{\pm}\pm i\xi _n)^{-j}$,
then$$
\multline
\int \bar{\hat\varphi }_k(\xi _n,\sigma ) h^+_{\xi _n}[\tilde
q(x',\xi ,\mu )\hat\varphi _l(\xi _n,\sigma )]\,\d \xi _n
=\int \bar{\hat\varphi }_k(\xi _n,\sigma ) \tilde
q(x',\xi ,\mu )\hat\varphi _l(\xi _n,\sigma )\,\d \xi _n\\
=\int h^-_{\xi _n}[\bar{\hat\varphi }_k(\xi _n,\sigma ) \tilde
q(x',\xi ,\mu )]\hat\varphi _l(\xi _n,\sigma )\,\d \xi _n.
\endmultline\tag B.5
$$
To see this, note that the integrand in all three expressions is
$O(\xi _n^{-2})$ for $|\xi _n|\to \infty $ in $\Bbb C$ and
meromorphic with no real poles, so
that the integral can be replaced by the integral over a large
contour in $\Bbb C_+$ (the ``plus-integral'') or the integral over a large
contour in $\Bbb C_-$.
The first equality holds since $\bar{\hat \varphi }_kh^- [\tilde
q\hat\varphi _l]$ is meromorphic with no poles in $\Bbb
C_+$ and is $O(\xi _n^{-2})$ for $|\xi _n|\to \infty $ in $\Bbb C$,
hence contributes with 0. The second equality holds since
$h^+[\bar{\hat \varphi }_k\tilde
q]\hat\varphi _l$ is meromorphic with no poles in $\Bbb C_-$, and is
$O(\xi _n^{-2})$ for $|\xi _n|\to \infty $ in $\Bbb C$, hence
contributes with 0.

The above formulas were used in \cite{GSc, Sect.\ 3 and 4} to show
that the trace expansions of terms where $G$ or $P_+$ is composed
with one of the parameter-dependent singular Green operators {\it do
not} contribute to the
residue coefficient $\tilde c_0'$ in (1.10). It was in fact shown that
they give expansions where the summation over $k$ as in (1.10) starts
with $k\ge 1$, so they do not contribute to the nonlocal coefficient
$\tilde c''_0$ either. This is important for our present study and
will therefore be formulated explicitly:

\proclaim{Proposition B.3} When $G$ is of class $0$ and order $\nu
\in\Bbb R$, or $P$ is of order $\nu \in\Bbb Z$, the operators
$GG^{(N)}_\lambda $,
$GG^{\pm}(Q^N_\lambda )$, $P_+G^{(N)}_\lambda $ and
$P_+G^{\pm}(Q^N_\lambda )$ have trace expansions of the form
$$
\sum_{j\ge
1}\!a_{j}(-\lambda)
^{\frac{n+\nu - j}{2}-N}
+\sum_{k\ge
1}(a'_k\log (-\lambda )+a''_k)(-\lambda )
^{-\frac k2-N}.\tag B.6
$$
\endproclaim
\demo{Proof} For the compositions with $G$ in front, this is
essentially the
content of \cite{GSc, Sect.\ 3}, see in particular (3.20) there. We
have replaced the index $l$ there by $j=l+1$ in the first sum,
$k=l+1$ in the second sum, to facilitate the comparison with (1.10);
recall also that $\mu =(-\lambda )^{\frac12}$.
The arguments there extend immediately to noninteger $\nu $.

For the compositions with $P_+$ in front, the statement is covered by the
calculations in \cite{GSc, Sect.\ 4}, see in particular (4.3) there
(which
contains neither nonlocal nor logarithmic terms, since $p'_{(l)}$ is
a differential operator) and (4.14) there.

In each case, the result is found by showing that in local
coordinates, the symbol of $\tr_n$ of the operator contains so many
negative powers of $\kappa ^{\pm}$ and $\kappa ^{\pm}+\sigma $ that
it is in $S^{\nu +1,-2N-1,0}\cap S^{\nu -2N,0,0}$, so that $d=-2N-1$
in (3.6).
\qed\enddemo

When $\nu \notin\Bbb Z$, the log-coefficients $a'_k$ vanish in (B.6),
since the degrees of the homogeneous symbols are noninteger.

The conclusions of the proposition hold also if the symbol $q$ of
$Q_\lambda $
is replaced by one of its $(x,\xi )$-derivatives, since they have a
similar structure.

The next lemma is used in  calculations of compositions of the
type $GQ^N_{\lambda ,+}$.

\proclaim{Lemma B.4}
Let
$$
s^{\pm}_{j,l,m}(y',\xi ',\mu )=\int \frac
{(\sigma - i\xi
_n)^{l}}{(\sigma + i\xi
_n)^{l+1}}\frac{(\sigma +
i\xi
_n)^{m}}{(\sigma -
i\xi
_n)^{m+1}}\frac1{(\kappa
^{\pm}\pm i\xi _n)^j}\,\d\xi _n.
\tag B.7
$$
One has for $l,m\in\Bbb Z$, $j\ge 1$:
$$\alignedat2
&\text{For $m<l$, }&&s^+_{j,l,m}=0,\quad
s^-_{j,l,m}=\sum_{ |l-m-1-m'|<j,\,  m'\ge
0}(\kappa^-) ^{1-j}b'_{jm'}\frac{(\sigma -\kappa ^-)^{m'}}{(\sigma +\kappa
^-
)^{m'+2}},\\
&\text{For $m=l$, }&&s^+_{j,l,l}=\frac1{(\kappa ^++\sigma )^j 2\sigma
},\quad s^-_{j,l,l}=\frac1{(\kappa ^-+\sigma )^j 2\sigma },\\
&\text{For $m>l$, }&&s^+_{j,l,m}=\sum_{ |m-l-1-m'|<j,\, m'\ge
0}(\kappa ^+) ^{1-j}b_{jm'}\frac{(\sigma -\kappa
^+)^{m'}}{(\sigma +\kappa ^+
)^{m'+2}},\quad s^-_{j,l,m}=0,
\endalignedat\tag B.8$$
with $b_{jm'}$ and $b'_{jm'}$ being $O(l^jm^j)$ for fixed $j$.

When $m\ne l$, the resulting symbols  are in $S^{0,0,-j-1}$
having $N$'th $\lambda $-derivatives in
$S^{0,0,-j-1-2N}$; for $m=l$, they are in  $S^{-1,0,-j}$,
with $N$'th $\lambda $-derivatives in $S^{1,0,-j-2N}$.

It follows (cf\. {\rm (A.12)}) that$$
\tr_n (\hat\varphi _l(\xi _n,\sigma )h^+(q^{\circ N}_{-2N}(x',\xi
',\eta _n,\mu
)\bar{\hat\varphi }_m(\eta _n,\sigma )))=\cases \alpha
^{(N)}(x',\xi ',\mu )&\text{ if }l=m,\\
s_{lm}(x',\xi ',\mu )&\text{ if }l\ne m,\endcases\tag B.9
$$
where $s_{lm}\in S^{1,0,-2N-1}$, the symbol seminorms being polynomially
bounded in
$l$ and $m$.
\endproclaim

\demo{Proof} The formulas in (B.8) were shown in \cite{GSc, Lemma
5.2}, and the
symbols are analyzed as in Lemma B.1.

Formula (B.9) is deduced from (B.7)--(B.8) using (B.5) and noting
that (B.7) contains
non-normalized Laguerre functions $\hat\varphi '_l$ and
$\bar{\hat\varphi }'_m$ so that we get an extra factor $2\sigma $
from the normalized Laguerre functions. The statements for $l\ne m$
follow straightforwardly from the descriptions in (B.8).
The statements for $l=m$ were proved in \cite{GSc, Prop\. 5.4};
they can be verified very simply by doing the calculation for $N=1$
and passing via $\lambda $-derivatives to the general case.
\qed
\enddemo

The important point in this lemma is that a symbol contributing to
$\tilde c'_0$ and $\tilde c''_0$ does appear, but with a special form
that allows further clarification (using the information (A.13)).

We summarize some results from \cite{GSc} on compositions of the form
$GQ^N_{\lambda,+} $ in the following statement:

\proclaim{Proposition B.5} Consider $GQ^N_{\lambda
,+}$ in a localization to $\Bbb R^n_+$.
Here  $$
\tr_n(GQ^N_{\lambda ,+})=\widetilde S_0+\widetilde S_1\text{ with }
\widetilde
S_0=\op'(\tr_ng(x',\xi
')\alpha ^{(N)}(x',\xi ',\mu )),
\tag B.10
$$
where
$\widetilde S_0$ and $\widetilde S_1$ are
$\psi $do's on $\Bbb R^{n-1}$ with symbols in
$S^{\nu ,-2N,0}\cap S^{\nu -2N,0,0}$ resp\.
\linebreak$S^{\nu +1,-2N-1,0}\cap S^{\nu -2N,0,0}$. The traces
$\Tr_{\Bbb R^{n-1}}(\widetilde S_0)$
 and $\Tr_{\Bbb R^{n-1}}(\widetilde S_1)$ have
expansions
$$
\sum_{j\ge
1}\!a_{j}(-\lambda)
^{\frac{n+\nu - j}{2}-N}
+\sum_{k\ge
0}(a'_k\log (-\lambda )+a''_k)(-\lambda )
^{-\frac k2-N},\tag B.11
$$
where the sum over $k$ starts at $k=1$ for $\widetilde S_1$, and
the value of $a'_0$ can be determined more precisely for $\widetilde S_0$.
\endproclaim

\demo{Proof}
This is proved in
\cite{GSc, Sect\. 5}, see in particular Prop.\ 5.3, 5.4,
and Sect\. 5.b there.\qed
\enddemo

We have in fact in view of Lemma A.3 that $\widetilde S_0\in
S^{\nu ,0,-2N}$, and also the information on $\widetilde S_1$
can be upgraded to $S^{\nu +1,0,-2N-1}$
by a closer analysis (as in Section 3 in this paper, using methods from
Proposition 3.8 to handle the case where $q$ depends on $x_n$).
\medskip

Let us introduce a notation for the Poisson and trace operators
on $\rnp$ (Laguerre operators) defined from Laguerre functions:
$$
\Phi _j=\opk(\hat{\varphi}_j(\xi_n,\sigma)),\quad \Phi ^*
_k=\opt(\bar{\hat{\varphi}}_k(\xi_n,\sigma)).\tag B.12
$$
Here $\Phi _j$ maps $L_2(\Bbb R^{n-1})$ continuously (in fact
isometrically) into $L_2(\rnp)$, and its adjoint is
$\operatorname{OPT}(\bar{\hat\varphi }_j)$. Moreover, because of
the orthonormality of the $\hat\varphi _j$,
$$
\Phi ^*_j\Phi
_k=\delta _{jk}I,\tag B.13
$$ where $I$ is the identity operator on functions on
$\Bbb R^{n-1}$.

\proclaim{Lemma B.6}
A singular Green operator $G$ on $\rnp$ of order $\nu$ and class 0
can be written in the form
$$
G=\sum_{j,k\in\Bbb N} \Phi _j C_{jk}
\Phi ^*_k
,\quad C_{jk}=\op'(c_{jk}(x',\xi')) , \tag B.14$$
with a rapidly decreasing sequence $(c_{jk})_{j,k\in\Bbb N}$  in
$S^{\nu}$.
\endproclaim

\demo{Proof} The symbol $g$ of $G$ has a Laguerre expansion:
$$
g(x',\xi',\xi_n,\eta_n) = \sum_{j,k\in\Bbb N}
d_{jk}(x',\xi'){\hat\varphi}_j(\xi_n,\sigma)
\bar{{\hat\varphi}}_k(\eta_n,\sigma),
\tag B.15$$
with $(d_{jk})_{j,k\in\Bbb N}$ rapidly decreasing in $S^{\nu}$,
i.e., the relevant symbol seminorms on \linebreak$\ang k^{N}\ang
j^{N'}d_{jk}$ are
bounded in $j,k$ for any $N,N'\in\Bbb N$. Then
$$
G=\sum_{k\in\Bbb N} \opk \big(\sum_{j\in\Bbb N}
d_{jk}(x',\xi'){\hat\varphi}_j(\xi_n,\sigma)\big)\circ
\opt\big(\bar{\hat\varphi}_k(\eta_n,\sigma)\big).
$$
Each Poisson operator
$\opk(\sum_{j\in\Bbb N}
d_{jk}(x',\xi'){\hat\varphi}_j(\xi_n,\sigma))$ has a symbol in $y'$-form,
$l_k(y',\xi)$. Since the $d_{jk}$ are rapidly decreasing in $S^{\nu}$,
the sequence $(l_k)$ is rapidly decreasing in the topology of Poisson
symbols
of order $\nu+\frac12$.
We now expand $l_k$ in a Laguerre series:
$$
l_k(y',\xi) = \sum_{j\in\Bbb N} c_{jk}(y',\xi')\hat\varphi_j(\xi_n,\sigma)
,
$$
and conclude from \cite{G2, Lemma 2.2.1} that $(c_{jk})$ is rapidly
decreasing in
$S^{\nu}$. Since
$$
\opk(l_k) = \sum_{j\in\Bbb N}
\opk(\hat\varphi_j(\xi_n,\sigma))\op'(c_{jk}(x',\xi')),
$$
we obtain $$
G=\sum_{j,k\in\Bbb
N}\opk\big(\hat\varphi_j(\xi_n,\sigma)\big)\op'\big(c_{jk}(x',\xi')\big)
\opt\big(\bar{\hat\varphi}_k(\eta_n,\sigma)\big),$$
which shows the assertion. \qed
\enddemo

The last part of the proof shows the following useful result:

\proclaim{Lemma B.7}
A Poisson operator $K$ of order $\nu+\frac12$ on $\overline{\Bbb R}^{n}_+$
can be written
$K=\sum_{j\ge 0}\Phi_jC_j$ with a rapidly decreasing sequence $(C_j)$ of
$\psi$do's
with symbols in $S^\nu$.
\endproclaim

\head{Acknowledgment} \endhead
Elmar Schrohe was partially supported by the European Research and Training Network
 `Geometric Analysis' (Contract HPRN-CT-1999-00118).

\enddocument